\documentclass[9pt]{article}
\usepackage{amsmath}
\usepackage{amsfonts}
\usepackage{bm}
\usepackage{mathrsfs}

\usepackage{amsthm}

\usepackage{leftidx}
\usepackage{graphicx}
\usepackage{hyperref}
\usepackage{xcolor}
\usepackage{graphicx, color, epstopdf}
\usepackage{cite}

\usepackage{subfigure}
\usepackage{empheq}

\usepackage{hyperref}
\newtheorem{theorem}{Theorem}[section]
\newtheorem{lemma}{Lemma}[section]
\newtheorem{proposition}{Proposition}[section]
\newtheorem{corollary}{Corollary}[section]

\newtheorem{remark}{Remark}

\renewenvironment{proof}[1][Proof]{\noindent\textit{#1. } }{\hfill$\square$}

\newcommand{\D}[1]{\mathcal{D}_{#1}}
\newcommand{\Ass}[1]{\textbf{\upshape A#1}}

\title{
A unconditionally energy dissipative, adaptive IMEX BDF2 scheme and its error estimates for Cahn-Hilliard equation on generalized SAV approach}

\author{
	Yifan Wei
 \thanks {School of Mathematics and Statistics, Wuhan University, Wuhan 430072, China (wei\_yi\_fan@whu.edu.cn)}
	\and Jiwei Zhang
	\thanks{School of Mathematics and Statistics, and Hubei Key Laboratory of Computational Science, Wuhan University, Wuhan 430072, China (jiweizhang@whu.edu.cn). }
	\and Chengchao Zhao
	\thanks{Beijing Computational Science Research Center, Beijing, 100193, P.R. China (cheng\_chaozhao@csrc.ac.cn).}
	\and Yanmin Zhao\thanks{Corresponding author. School of Science, Xuchang University, Xuchang, 461000, China. (zhaoym@lsec.cc.ac.cn)}
	}
\date{}

\begin{document}

\maketitle

\begin{abstract}
 An adaptive implicit-explicit (IMEX) BDF2 scheme is investigated on generalized SAV approach for the Cahn-Hilliard equation by combining with Fourier spectral method in space.
It is proved that the modified energy dissipation law is unconditionally preserved at discrete levels.
Under a mild ratio restriction,  i.e., \Ass{1}: $0<r_k:=\tau_k/\tau_{k-1}< r_{\max}\approx 4.8645$, we establish a rigorous error estimate in $H^1$-norm and achieve optimal second-order accuracy in time.
The proof involves the tools of discrete orthogonal convolution (DOC) kernels and inequality zoom.
It is worth noting that the presented adaptive time-step scheme only requires solving one linear system with constant coefficients at each time step.
In our analysis, the first-consistent BDF1 for the first step does not bring the order reduction in $H^1$-norm.
The $H^1$ bound of numerical solution under periodic boundary conditions can be derived  without any restriction (such as zero mean of the initial data).
Finally, numerical examples are provided to verify our theoretical analysis and the algorithm efficiency.
\end{abstract}


\section{Introduction}
The Cahn-Hilliard (C-H) equation, proposed in \cite{cahn_free_nodate} to model the process of phase separation in binary alloys,  has been widely used in diblock copolymer \cite{zhang_new_2020}, image inpainting \cite{brkic_image_2020}, tumor growth simulation \cite{agosti_cahn-hilliard-type_2017} and topology optimization \cite{bartels_cahnhilliard_2021}.
In this paper, we consider the computation of  the following C-H equation 
\begin{equation}\label{CHeq}
 \partial_t\Phi=\Delta \mu \quad \text{with} \quad \mu=-\Delta \Phi+ \frac{1}{\varepsilon^2}(\Phi^3-\Phi),\quad (\bm{x},t)\in\Omega\times(0,T],
\end{equation}
with periodic boundary and the initial condition $\Phi(\bm{x} , 0) = \Phi^0(\bm{x})$. Here $\Omega \in \mathbb{R}^d$ ($d= 2,3$), $\mu$ is the chemical potential and the positive parameter $\varepsilon$ measures the width of the diffuse interface layer. The system \eqref{CHeq}  can be viewed as the $H^{-1}$ gradient flow of the Ginzburg-Landau type energy functional
\begin{equation}\label{GL}
 E[\Phi]=\int_\Omega\left( \frac{1}{2} |\nabla \Phi|^2+\frac{(\Phi^2-1)^2}{4\varepsilon^2}\right) \textrm{d}\bm{x},
\end{equation}
which holds the following energy dissipation law
\begin{align}\label{EnDissipation}
\frac{\textrm{d}}{\textrm{d} t}E[\Phi(t)]=\int_\Omega \frac{\delta E}{\delta \Phi}\partial_t\Phi\textrm{d}\bm{x}=-\int_\Omega |\nabla \mu|\textrm{d}\bm{x}.
\end{align}

In any numerical methods for solving C-H model, a key consideration is preserving the dissipativity, which can capture the long-time statistical properties under the approximation.
Many efforts have been made to construct efficient numerical schemes that satisfy energy dissipation at discrete levels, possibly in some modified form \cite{liao_adaptive_2020}.
These includes, but not limited to, the convex splitting method \cite{eyre_unconditionally_1998,chen_second_2019, cheng_energy_2019}, linear stabilization approach \cite{MR2322435,MR3109557}, the invariant energy quadratization (IEQ) \cite{yang_numerical_2017} and the scalar auxiliary variable (SAV) \cite{shen_scalar_2018}.
In particular, a generalized SAV approach (gSAV) \cite{huang_highly_2020,MR4383075} offers essential improvements over the original SAV approach.
Specifically, the gSAV only requires solving one linear system with constant coefficients at each time step, which means the computational cost is about half of the original SAV.
Meanwhile, it applies to more general gradient flows and even to general dissipative systems.
But, the convergence analysis in \cite{huang_highly_2020} is only for temporal uniform mesh.

Noticing that the multi-scale feature of the C-H equation, it is worthwhile to consider the adaptive mesh in long-time simulations.
In fact, there are fast time scales in the dynamics of many phase field models such as C-H equation. Very fine time step size is needed to capture the underlying physical or biological phase transition phenomena.
Meanwhile, there are long periods of slow phase transitions  before equilibrium is reached. In this situation, relatively large step sizes can maintain the same accuracy.
Therefore, it is highly desired to develop adaptive time-stepping techniques which improves efficiency without sacrificing accuracy \cite{qiao_adaptive_2011,MR3072975}.
Any convincing adaptive strategy needs theoretical support and guidance.
On numerical analysis of adaptive multi-step methods, such as the BDF2 scheme, are quite difficult compared with one-step methods.
In fact, Becker \cite{becker_second_1998} presents the bounds of stability and second-order convergence under adjacent ratio restriction $r_k\leq 1.8685$ and the boundedness of $\sum_{k=3}^n\max\{0,r_{k-1}-r_{k+1}\}$.
To circumvent this strict restriction, a new developed Gr\"onwall inequality was employed in \cite{chen_second_2019} and the convergence analysis was derived under the restriction $r_k\leq 1.534$.
Recently, the restriction on the adjacent time step ratio have been extend to $r_k\leq3.561$ in \cite{liao_analysis_2021-1} and $r_k\leq4.8645$ in \cite{zhang_sharp_2020} for linear parabolic equations by using the techniques of DOC/DCC kernels.
 Furthermore, 
the BDF2 method is A-stable and L-stable such that it is widely used for solving the stiff problems \cite{MR3997384,fan}.

Many works have been carried out on energy stable variable-step BDF2 schemes for phase field models, most of them employ implicit schemes (fully implicit or partially implicit).
For instance, a fully implicit variable-step BDF2 scheme is studied in \cite{fan} to achieve the robust second-order convergence for the phase field crystal model.
Relying on convex splitting and stabilization technique, unconditional energy stable variable-step BDF2 schemes are constructed and analyzed for C-H equation \cite{chen_second_2019,MR4448198}.
For the extended Fisher-Kolmogorov equation, a variable-step BDF2 scheme is established in \cite{MR4401877} using the convex splitting strategy.
Note that the fully/partially implicit schemes generally require more computational costs as they need to solve a resulting nonlinear system at each time step.
In contrast, linear schemes \cite{MR2948706_Chen,MR3199497_Chen2,MR3564340_Yang}  only need to solve linear systems at each time step. Thus, an adaptive IMEX BDF2 scheme  will be more efficient in practical simulations.
So far,  few work has been done to study linear schemes on temporal adaptive mesh with energy dissipation laws.
Note that the gSAV method has many advantages over traditional SAV  \cite{huang_highly_2020} such as reduction in computational costs.
It is worthwhile to study an adaptive IMEX  BDF2 scheme with gSAV and its convergence analyses under a mild ratio restriction as
\begin{center}
\Ass{1} : \quad $0< r_k \leq r_{\max} (\approx 4.864)-\delta$,
\end{center}
where $r_{\max}$ is the real root of $x^3=(2x+1)^2$ and $\delta$ is any given small constant such as $\delta = 0.01$.

In this paper, we construct and analyze an adaptive IMEX BDF2 scheme on gSAV approach for C-H equations which enjoys the advantage of solving one linear system with constant coefficients at each time step.
The main contribution of this paper has two aspects:
(i) we rigorously prove the unconditionally energy dissipation property  of the proposed variable time-step BDF2 scheme;
(ii) we establish  the  corresponding convergence analysis in $H^1$-norm under the mild ratio restriction \Ass{1} while  the second-order optimal error estimation is achieved in time.
The proof involves the tools of DOC kernels and its generalized properties.
One of the difficulties of error analysis is to prove rigorously that first-consistent BDF1 for the first step does not result in the order reduction in $H^1$-norm (reference \cite{MR4206656} for detail about order reduction).
This can be obtained by inequality zoom in the proof and some delicate error analysis on truncation error (see Remark \ref{remark2}).
Besides, one of the main defects in existing theory is  that the $H^1$ bound is proved to require the initial data having mean zero \cite{MR4379976,shen_convergence_2018} or adding an extra term $\frac{\lambda}{2} u^2$  \cite{shen_convergence_2018}.
In this work, these extra restrictions are circumvented by a novel technique (see Remark \ref{remark1}).


The rest of the paper is organized as follows.
In section 2, a fully discrete IMEX BDF2 scheme with variable-time step  is presented.
In section 3, the modified energy dissipation law of the proposed scheme and the $H^1$ bounded of $\phi^n$ are considered.
On this basis, the $H^2$ bounded of $\phi^n$ under the regularity of $\Phi^0\in H^4$ is proved.
The $H^1$-norm convergence analysis is established by mathematical induction in section 4. In the last section, numerical experiments are provided to demonstrate the convergence, energy dissipation properties and efficiency of adaptive strategy, respectively.

\section{Preliminaries}
The variable time mesh is given as  $0=t_0<t_1<\cdots <t_{K}=T$. Denote  the time step  by $\tau_k:= t_k-t_{k-1}\ (k\geq 1)$, the maximum time step by $\tau:=\max_{1\leq k\leq K}\tau_k$ and the adjacent time-step ratio by $r_k:=\tau_k/\tau_{k-1}\ (k\geq 2)$ with $r_1\equiv 0$.
Define the difference operator by $\nabla_\tau u^k:= u^k-u^{k-1}$. Then the variable time-step BDF1 and BDF2 formulas are respectively defined by
\begin{equation*}
 \mathcal{D}_1 u^k= \frac{1}{\tau_k} \nabla_\tau u^k, \quad
 \mathcal{D}_2 u^k  = \frac{1+2r_k}{\tau_k(1+r_k)}\nabla_\tau u^k - \frac{r_k^2}{\tau_k(1+r_k)}\nabla_\tau u^{k-1}.
\end{equation*}
Since BDF2 is a two-step method requiring two starting values, we use the  BDF1 to compute the first-step value. Introducing the discrete convolution kernels $b^{(n)}_{n-k}$ as follows
\begin{equation*}
b^{(n)}_{0}:=\frac{1+2r_n}{\tau_n(1+r_n)},\  b^{(n)}_{1}:=-\frac{r_n^2}{\tau_n(1+r_n)}\quad  \text{and} \quad  b^{(n)}_{j}:=0 \quad \text{for}\ \  j\geq 2,
\end{equation*}
the BDF2 started by the BDF1 can be written as a unified convolution form
\begin{equation}\label{OpeaD2}
\mathcal{D}_2 u^n:=\sum\limits^n_{k=1}b^{(n)}_{n-k}\nabla_\tau u^k\ \  \text{for} \ \ n\geq 1.
\end{equation}
We further define the extrapolation operator $B$ as
\begin{equation}\label{OpeaB}
Bu^{n-1}=(1+r_n)u^{n-1}-r_nu^{n-2}\quad\forall n\geq 2 \quad \text{and} \quad Bu^0= u^0.
\end{equation}


 We adopt the Fourier spectral method for the spatial discretization of \eqref{CHeq}.
For simplicity,  we take the two-dimensional domain $\Omega=(0, L)^2$ as an example, which is  partitioned by uniform mesh size $h=L/N$. Here $N$ denoting the numbers of Fourier modes in each direction.  Define Fourier approximation space by
 $$S_N= \text{span}\{e^{\textrm{i}\xi_kx}e^{\textrm{i}\eta_ly}:-\frac{N}{2} \leq k \leq \frac{N}{2}-1, -\frac{N}{2} \leq l \leq \frac{N}{2}-1 \},$$
where $\textrm{i}=\sqrt{-1}, \xi_k=2\pi k/L$ and $\eta_l=2\pi l/L$.
The function $u(x,y)\in L^2(\Omega)$ can be approximated by
\begin{equation*}
u(x,y)\approx u_N(x,y)=\sum\limits_{k=-N/2}^{N/2-1}\sum\limits_{l=-N/2}^{N/2-1} \hat{u}_{k,l}e^{\textrm{i}\xi_k x}e^{\textrm{i}\eta_l y},
\end{equation*}
where the Fourier coefficients $\hat{u}_{k,l}$ are given as
$\hat{u}_{k,l}=\frac{1}{|\Omega|}\int_\Omega ue^{-\textrm{i}(\xi_k x+\eta_l y)} \mathrm{d}x\mathrm{d}y.$

{Denote $(\cdot,\cdot)$ and $\|\cdot\|$ by  the inner product and  norm  in $L^2(\Omega)$ respectively.}
The $L^2$-orthogonal projection operator $P_N:L^2(\Omega)\rightarrow S_N$ is defined by
\begin{equation}\label{OpeaP}
(P_N u-u,v_N)=0, \quad \forall v_N\in S_N, u\in L^2(\Omega).
\end{equation}
For the construction of SAV scheme, we now introduce a modified energy $\Gamma (t) =E(\Phi)+1$.
Then the equation \eqref{CHeq} with the energy dissipation law \eqref{EnDissipation} can be rewritten as a equivalent form:
\begin{align}
&(\partial_t\Phi,v)-(\Delta \mu,v)=0, \quad \forall v\in L^2(\Omega), \label{SAV1}\\
&\mu=\frac{\delta E}{\delta \Phi}=-\Delta\Phi+f(\Phi),\\
&\frac{\textrm{d} \Gamma}{\textrm{d} t}=-\frac{\Gamma}{E+1}\|\nabla\mu\|^2,\label{SAV2}
\end{align}
where $f(\Phi)=\frac{1}{\varepsilon^2}(\Phi^3-\Phi).$
Noting the periodic boundary conditions, it follows that \eqref{CHeq} satisfies mass conservation (i.e., $\int_\Omega \partial_t \Phi \textrm{d}\bm{x} \equiv0$) by choosing $v=1$ in \eqref{SAV1} .

A linear second-order scheme on generalized SAV approach  for problem \eqref{SAV1}-\eqref{SAV2} is constructed with Fourier spectral method in space as follows.
Given $\phi^{n-1},\phi^{n-2},\bar\phi^{n-1},\bar\phi^{n-2}\in S_N,\gamma^{n-1}\in \mathbb{R}$, we compute $\bar\phi^n, \gamma^n, \xi^n, \eta^n, \phi^n$  in sequence by
\begin{subequations}\label{scheme}
\begin{empheq}[left=\empheqlbrace]{align}
&(\D 2 \bar\phi^n,v_N)+(\Delta^2\bar\phi^n, v_N)-(\Delta f(B\phi^{n-1}),v_N)=0, \ \forall v_N\in S_N,\label{scheme1}\\
&\frac{\gamma^n-\gamma^{n-1}}{\tau_n}=-\frac{\gamma^n}{E(\bar\phi^n)+1}\|-\nabla\Delta\bar\phi^n+\nabla f(B\phi^{n-1})\|^2,\label{scheme2}\\
&\xi^n=\frac{\gamma^n}{E(\bar\phi^n)+1},\label{scheme3}\\
&\phi^n=\eta^n\bar\phi^n \quad
\text{where} \quad\eta^n=1-(1-\xi^n)^2 = \xi^n(2-\xi^n),\label{scheme4}
\end{empheq}
\end{subequations}
where the initial values are given by $\bar\phi^0=\phi^0=P_N \Phi^0$ and $\gamma^0=E(\phi^0)+1$. The operators $\D 2$ and $B$ are defined in \eqref{OpeaD2} and \eqref{OpeaB}, respectively.
Choosing $v_N=1$ in \eqref{scheme1}, we can indicate that the fully discrete scheme \eqref{scheme} satisfies mass conservation in the sense of
\begin{equation}
(\bar\phi^0,1)=(\bar\phi^n,1), \quad \forall n\geq 1.\label{Conser}
\end{equation}
Through this paper, $C, C_\Omega$ denote positive constant, and are not necessarily the same at different occurrences, but independent of the parameters and functions involved.

\section{Energy dissipation law and numerical stability}
A modified energy dissipation property at the discrete levels is considered in this section. Meanwhile,  the $H^1$ and $ H^2$ bounds of $\phi^n$ is rigorously proved.
\subsection{Discrete energy dissipation law}
\begin{lemma}\label{Lem_3.1} It holds that
\begin{equation}
\frac{\alpha x}{x^4+\alpha} < \alpha^{1/4}, \quad  \forall \alpha>0,\, x\geq0. \label{eq_L3}
\end{equation}
\end{lemma}
\begin{proof}
It is trivial to see that \eqref{eq_L3} holds when $x=0$. One only needs to consider the case of $x>0$.
Denote $g(x)=x^3+\alpha/x$. The direct calculation shows that
\begin{equation*}
g'(x)=3x^2-\frac{\alpha}{x^2},\quad
g''(x)=6x+\frac{2\alpha}{x^3}> 0 \text{ and }
 g'((\frac{\alpha}{3})^{1/4})=0,
\end{equation*}
which implies the  minimum
$\min\limits_{x>0}g(x)=g((\frac{\alpha}{3})^{1/4})$.
Furthermore, it holds
\begin{equation*}\frac{\alpha x}{x^4+\alpha}=\frac{\alpha}{g(x)}\leq \frac{\alpha}{\min\limits_{x>0}g(x)}
=\frac{\alpha^{1/4}}{3^{1/4}+3^{-3/4}} < \alpha^{1/4}.
\end{equation*}
The proof is completed.
\end{proof}

Noticing $L^4(\Omega)\hookrightarrow L^2(\Omega)$, Sobolev's embedding theorem implies there exists a constant $\bar{c}_\Omega$ such that
\begin{equation}\label{Sob_embed}
  \|u\|_{L^2}\leq \bar{c}_\Omega\|u\|_{L^4},\quad \forall u\in L^4(\Omega).
\end{equation}

We now consider the energy  dissipation with respect to  modified discrete energy $\gamma^n$ and  the $H^1$ bound of numerical solutions $\phi^n$ as follows.
\begin{theorem}\label{TH3.1}
Assume $\gamma^n, \phi^n$ are the solutions of \eqref{scheme},
then $\gamma^n>0, \xi^n>0$ and the scheme \eqref{scheme}
preserves the energy dissipation law in the sense that
\begin{equation}\label{eq_T3}
\gamma^n-\gamma^{n-1}\leq -\tau_n\xi^n\|-\nabla\Delta\bar\phi^n+\nabla f(B\phi^{n-1})\|^2\leq 0.
\end{equation}
Furthermore, it holds
$$\|\phi^n\|_{H^1}\leq \mathcal{M},$$
where $\mathcal{M}=2\gamma^0(\gamma^0+2)\sqrt{1+\bar{c}_\Omega^2\sqrt{2\varepsilon^2+\|1\|^2}}$ and the constant $\bar{c}_\Omega$ is given in \eqref{Sob_embed}.
\end{theorem}

\begin{proof}
 Eq. \eqref{scheme2} can be rewritten as
\begin{equation}\label{eq_T3.1}
\gamma^n(1+\frac{\tau_n\|-\nabla\Delta\bar\phi^n+\nabla f(B\phi^{n-1})\|^2}{E(\bar\phi^n)+1})
=\gamma^{n-1}.
\end{equation}
In view of $\gamma^0=E(\phi^0)+1>0$, it directly follows from \eqref{eq_T3.1} that
$$0<\gamma^n\leq\gamma^{n-1}\leq \cdots \leq \gamma^1\leq \gamma^0, \quad \forall n\geq 2.$$
Then, by use of \eqref{scheme2}-\eqref{scheme3}, one has  \eqref{eq_T3} and
\begin{equation}
  0< \xi^n = \frac{\gamma^n}{E(\bar\phi^n)+1} \leq \frac{\gamma^0}{E(\bar\phi^n)+1}\leq \gamma^0.\label{eq_T3.2}
\end{equation}

We now consider the $H^1$ bound for $\phi^n$.
From \eqref{eq_T3.2} and the definition of $\eta^n$, we have
\begin{equation}\label{eq_T3.3}
|\eta^n|=|\xi^n(2-\xi^n)| \leq  \frac{\gamma^0(|\xi^n|+2)}{E(\bar\phi^n)+1} \leq  \frac{\gamma^0(\gamma^0+2)}{E(\bar\phi^n)+1}.
\end{equation}
Considering \eqref{scheme4} and \eqref{eq_T3.3}, we have
\begin{equation}
\| \phi^n\|_{H^1}^2 \leq|\eta^n|^2\|\bar\phi^n\|_{H^1}^2
\leq [\gamma^0(\gamma^0+2)]^2\frac{\|\bar\phi^n\|^2+\|\nabla\bar\phi^n\|^2}{(E(\bar\phi^n)+1)^2}:=[\gamma^0(\gamma^0+2)]^2(M_1^2+M_2^2).\label{eq_T3.4}
\end{equation}
Noticing that $E(\bar\phi^n)=\frac{1}{2}\|\nabla\bar\phi^n\|^2+\frac{1}{4\varepsilon^2}\|(\bar\phi^n)^2-1\|^2$ and \eqref{Sob_embed}, we arrive at
\begin{equation}
M_1=\frac{\|\nabla\bar\phi^n\|}{E(\bar\phi^n)+1}\leq\frac{2\|\nabla\bar\phi^n\|}{\|\nabla\bar\phi^n\|^2+2}\leq 2,\quad
M_2=\frac{\| \bar\phi^n\|}{E(\bar\phi^n)+1}\leq \frac{4\varepsilon^{2}\bar{c}_\Omega\|\bar\phi^n\|_{L^4}}{\|(\bar\phi^n)^2-1\|^2+4\varepsilon^2}.\label{eq_T3.5}
\end{equation}
To deduce $\|\phi^n\|_{H^1}\leq \mathcal{M}$, we need to prove $M_1^2+M_2^2\leq 4(1+\bar{c}_\Omega^2\sqrt{2\varepsilon^2+\|1\|^2})$ in \eqref{eq_T3.4}. To this end,  from \eqref{eq_T3.5} we only need to prove
\begin{equation}\label{eq_T3_Inequal}
\frac{2\varepsilon^{2}\|\bar\phi^n\|_{L^4}}{\|(\bar\phi^n)^2-1\|^2+4\varepsilon^2}\leq  (2\varepsilon^2+\|1\|^2)^{\frac{1}{4}}.
\end{equation}
Along with $(a+b)^2\leq (1+\epsilon)a^2+(1+\frac{1}{\epsilon})b^2$ for any $\epsilon>0$,
we have
\begin{align}
\|\bar\phi^n\|^4_{L^4}= \|(\bar\phi^n)^2\|^2=\|((\bar\phi^n)^2-1)+1\|^2&\leq (1+\frac{\|1\|^2}{2\varepsilon^2})\|(\bar\phi^n)^2-1\|^2+(1+\frac{2\varepsilon^2}{\|1\|^2})\|1\|^2\nonumber\\
& = \frac{2\varepsilon^2+ \|1\|^2}{2\varepsilon^2}\|(\bar\phi^n)^2-1\|^2+\|1\|^2+ 2\varepsilon^2.\label{eq_T3.6}
\end{align}
It follows from \eqref{eq_T3.6} that the denominator of \eqref{eq_T3_Inequal} can be estimated by
\begin{align}\label{eq_T3.7}
\|(\bar\phi^n)^2-1\|^2+4\varepsilon^2&\geq 
\frac{2\varepsilon^2}{2\varepsilon^2+\|1\|^2}(\|\bar\phi^n\|^4_{L^4}+2\varepsilon^2+\|1\|^2).
\end{align}
Inserting \eqref{eq_T3.7} into \eqref{eq_T3_Inequal}, we have
\begin{equation*}
\frac{2\varepsilon^{2}\|\bar\phi^n\|_{L^4}}{\|(\bar\phi^n)^2-1\|^2+4\varepsilon^2}
\leq \frac{\|\bar\phi^n\|_{L^4}(2\varepsilon^2+\|1\|^2)}{\|\bar\phi^n\|^4_{L^4}+2\varepsilon^2+\|1\|^2}
< (2\varepsilon^2+\|1\|^2)^{\frac{1}{4}},
\end{equation*}
where Lemma \ref{Lem_3.1} is used by taking $x=\|\bar\phi^n\|_{L^4}$ and $\alpha = 2\varepsilon^2+\|1\|^2$.
The proof is completed.
\end{proof}


\begin{remark}\label{remark1}
We point out that additional restrictions are required in the existing SAV method to derive the $H^1$ bound by the modified discrete energy $\gamma$.
For instance, the initial values are restricted to zero mean (i.e., initial total mass satisfying $\int_\Omega \Phi^0 \rm{d}x=0$) in \cite{huang_stability_2021,MR4379976}. In \cite{shen_convergence_2018}, additional term $\frac{\lambda}{2}u^2$ is added for the Ginzburg-Landau type energy functional \eqref{GL}, which ensures $\frac{\lambda}{2}\|u\|^2+\frac{1}{2}\|\nabla u\|^2$ is a norm in $H^1$.
{In this paper, the proof of Theorem \ref{TH3.1} does not requires the  additional restrictions in \cite{shen_convergence_2018,huang_stability_2021,MR4379976}.}

\end{remark}

\subsection{\texorpdfstring{$H^2$}. bound of \texorpdfstring{$\phi^n$}. }
We now consider the $H^2$ bound by introducing the  discrete orthogonal convolution (DOC) kernel as
\begin{equation}\label{DOC}
 \sum_{j=k}^n \theta_{n-j}^{(n)} b_{j-k}^{(j)} = \delta_{nk}, \quad \forall  1\leq k \leq n,  \; 1 \leq n \leq K,
\end{equation}
where $\delta_{nk} = 1$ if $n=k$ and $ \delta_{nk} = 0$ if $n\neq k$.
 According to  definition \eqref{DOC},  we have
 \begin{equation} \label{Theta_D2}
\sum_{j=1}^n\theta_{n-j}^{(n)}\mathcal{D}_2 u^j= \sum_{l=1}^n\nabla_\tau u^l \sum_{j=l}^n\theta_{n-j}^{(n)}b_{j-l}^{(j)}  =  u^n - u^{n-1}, \quad \forall 1\leq n \leq K.
\end{equation}

\begin{lemma}[\!\!\cite{mayuheng2}]\label{Lem_inequal}
If $\mathbf{A1}$ holds, then for any real sequences $\{w_k\}_{k=1}^n$, it has
\begin{equation*}
2\sum\limits_{k=1}^n w_k\sum\limits_{j=1}^k\theta^{(k)}_{k-j} w_j
\geq \sum\limits_{k=1}^n \frac{\delta}{20}\frac{(\sum^n_{s=k}\theta^{(s)}_{s-k} w_s)^2}{\tau_k}\geq c_\delta\sum\limits_{k=1}^n\tau_k w_k^2 \geq 0,
\text{ for }n\geq 1,
\end{equation*}
where $c_\delta$ is a constant dependent on $\delta$.
\end{lemma}

\begin{lemma} [\!\!\cite{shen_convergence_2018, 1900Infinite}]\label{Lem_3.3}
Assume that $\|u\|_{H^1}\leq  M$ and
\begin{align*}
|g'(x)|\leq C(|x|^{p_1}+1), p_1>0\text{ arbitrary if }d=1,2;\\
|g''(x)|\leq C(|x|^{p_2}+1), p_2>0\text{ arbitrary if }d=1,2.
\end{align*}
Then, for any $u\in H^4$, there exist $0\leq \sigma<1$ and a constant $C(M)$ such that
$$\|\Delta g(u)\|^2\leq C(M)(1+\|\Delta^2 u\|^{2\sigma}).$$
Furthermore, for any $\epsilon > 0$, there exists a constant $C(\epsilon, M)$ depending on $\epsilon$ such that
\begin{equation}
\|\Delta g(u)\|^2\leq \epsilon\|\Delta^2 u\|^2+C(\epsilon,M).
\label{Lem_3.3_Y}
\end{equation}
\end{lemma}
Then we give an $H^2$ bound for \eqref{scheme} as follows.
\begin{theorem}\label{TH3.2_bound}
Assume  $\bar{\phi}^n$ and $\phi^n$ are the solutions of \eqref{scheme} and $\mathbf{A1}$ holds.
For any $\epsilon > 0$, it holds
\begin{equation}
\|\Delta\bar\phi^n\|^2+\big{(}\frac{c_\delta}{2}-\frac{40\epsilon c_1}{\delta}\big{)}\sum\limits_{k=1}^{n-1}\tau_k\| \Delta^2\bar \phi^k\|^2\leq \|\Delta\bar{\phi}^0\|^2+ \frac{80(1+r_{\max})^3\epsilon\tau_1}{\delta}\|\Delta^2\phi^0\|^2+\frac{C(\epsilon,M)t_n}{\delta},\label{eq_TH32_1}
\end{equation}
where
\begin{equation}
c_1=2(1+r_{\max})^3r_{\max}(\gamma^0)^2(2+\gamma^0)^2.\label{c1}
\end{equation}
\end{theorem}
\begin{proof}
Choosing $v_N=2\Delta^2\bar\phi^k$ in \eqref{scheme1}, we have
\begin{equation}\label{eq_TH32_2}
2(\D 2 \bar\phi^j,\Delta^2\bar\phi^k)+2(\Delta^2\bar\phi^j,\Delta^2\bar\phi^k)-2(\Delta f(B\phi^{j-1}),\Delta^2\bar\phi^k)=0.
\end{equation}
Multiplying \eqref{eq_TH32_2} by $\theta^{(k)}_{k-j}$ and sum over $j$ from 1 to $k$, we then sum over $k$ from 1 to $n$ to get
\begin{equation}\label{eq_TH32_3}
2\sum\limits_{k=1}^n(\sum^k_{j=1}\theta^{(k)}_{k-j}\D 2 \Delta\bar\phi^j,\Delta\bar\phi^k)
+2\sum\limits_{k=1}^n\sum\limits_{j=1}^k\theta^{(k)}_{k-j}(\Delta^2\bar \phi^j,\Delta^2\bar \phi^k)
= 2\sum\limits_{k=1}^n\sum\limits_{j=1}^k\theta^{(k)}_{k-j}(\Delta f(B\phi^{j-1})
,\Delta^2\bar \phi^k).
\end{equation}
Applying \eqref{Theta_D2} and $2a(a-b)\geq a^2-b^2$,
the first term on the left-hand side of \eqref{eq_TH32_3} can be estimated as
\begin{equation}\label{eq_TH32_4}
2\sum\limits_{k=1}^n(\sum^k_{j=1}\theta^{(k)}_{k-j}\D 2 \Delta\bar\phi^j,\Delta\bar\phi^k)\geq\|\Delta\bar\phi^n\|^2-\|\Delta\bar\phi^{0}\|^2.
\end{equation}
It follows from Lemma \ref{Lem_inequal} that
\begin{equation}
2\sum\limits_{k=1}^n\sum\limits_{j=1}^k\theta^{(k)}_{k-j}(\Delta^2\bar \phi^j,\Delta^2\bar \phi^k)
\geq \sum\limits_{k=1}^n \frac{\delta}{20}\frac{\|\sum^n_{s=k}\theta^{(s)}_{s-k} \Delta^2\bar \phi^s\|^2}{\tau_k}.
\end{equation}
Exchanging the summation order and applying Young's inequality, we have
\begin{align}\label{eq_TH32_5}
2\sum\limits_{k=1}^n\sum\limits_{j=1}^k\theta^{(k)}_{k-j}(\Delta f(B\phi^{j-1})
,\Delta^2\bar \phi^k)
=2\sum\limits_{j=1}^n(\Delta f(B\phi^{j-1})
,\sum\limits_{k=j}^n\theta^{(k)}_{k-j}\Delta^2\bar \phi^k)\nonumber\\
\leq 40\delta^{-1}\sum\limits_{j=1}^n\tau_j\|\Delta f(B\phi^{j-1})\|^2
+\delta\sum\limits_{j=1}^n\frac{\|\sum_{k=j}^n\theta^{(k)}_{k-j}\Delta^2\bar \phi^j\|^2}{40\tau_j}.
\end{align}
Inserting \eqref{eq_TH32_4}-\eqref{eq_TH32_5} into \eqref{eq_TH32_3}, we obtain
\begin{equation}\label{eq_TH32_6}
\|\Delta\bar\phi^n\|^2
+\delta\sum\limits_{j=1}^n\frac{\|\sum_{k=j}^n\theta^{(k)}_{k-j}\Delta^2\bar \phi^j\|^2}{40\tau_j}
\leq \|\Delta\bar\phi^{0}\|^2+ 40\delta^{-1}\sum\limits_{k=1}^n\tau_k\|\Delta f(B\phi^{k-1})\|^2.
\end{equation}
It follows from Theorem \ref{TH3.1} to have $\|\phi^n\|_{H^1}\leq  M$.
Applying  \eqref{Lem_3.3_Y} of Lemma \ref{Lem_3.3} yields
\begin{equation}\label{eq_TH32_7}
\sum\limits_{k=1}^n\tau_k\|\Delta f(B\phi^{k-1})\|^2\leq \epsilon \sum\limits_{k=1}^n\tau_k\|\Delta^2 B\phi^{k-1}\|^2+ C(\epsilon,M)\sum\limits_{k=1}^n\tau_k .
\end{equation}
Noticing that
\begin{equation*}
\|\Delta^2 B \phi^{k-1}\|^2=\|(1+r_k)\Delta^2\phi^{k-1}-r_k \Delta^2\phi^{k-2}\|^2\leq 2(1+r_{\max})^2(\|\Delta^2\phi^{k-1}\|^2+\|\Delta^2\phi^{k-2}\|^2), \quad k\geq 2,
\end{equation*}
and $\|\Delta^2B \phi^0\|^2=\|\Delta^2\phi^0\|^2$, we have
\begin{align}\label{eq_TH32_8}
\sum\limits_{k=1}^n\tau_k \|\Delta^2 B\phi^{k-1}\|^2
&\leq 2(1+r_{\max})^3\sum\limits_{k=1}^n\tau_k\|\Delta^2\phi^{k-1}\|^2\nonumber\\
&\leq 2(1+r_{\max})^3r_{\max}\sum\limits_{k=1}^{n-1}\tau_{k}\|\Delta^2\phi^{k}\|^2
+2(1+r_{\max})^3\tau_1\|\Delta^2\phi^{0}\|^2.
\end{align}
It follows from \eqref{eq_T3.3} that $|\eta^n|\leq \gamma^0(2+\gamma^0)$. Then, it is a consequence of \eqref{scheme4} that
\begin{equation}\label{eq_TH32_9}
\sum\limits_{k=1}^{n-1}\tau_{k}\|\Delta^2 \phi^{k}\|^2
\leq (\gamma^0)^2(2+\gamma^0)^2\sum\limits_{k=1}^{n-1}\tau_k\| \Delta^2 \bar \phi^k\|^2.
\end{equation}
Combining \eqref{eq_TH32_7}-\eqref{eq_TH32_9}, we arrive at
\begin{equation}\label{eq_TH32_10}
\sum\limits_{k=1}^n\tau_k\|\Delta f(B\phi^{k-1})\|^2
\leq  \epsilon c_1\sum\limits_{k=1}^{n-1}\tau_k\| \Delta^2\bar \phi^k\|^2+2(1+r_{\max})^3\epsilon\tau_1\|\Delta^2\phi^{0}\|^2+ C(\epsilon,M)t_n,
\end{equation}
where $c_1$ is defined in \eqref{c1}.
By Lemma \ref{Lem_inequal}, we have
\begin{equation}\label{eq_TH32_11}
\delta\sum\limits_{j=1}^n\frac{\|\sum_{k=j}^n\theta^{(k)}_{k-j}\Delta^2\bar \phi^j\|^2}{40\tau_j}
\geq \frac{c_\delta}{2}\sum\limits_{k=1}^{n-1}\tau_k\| \Delta^2\bar \phi^k\|^2.
\end{equation}
Inserting \eqref{eq_TH32_10}, \eqref{eq_TH32_11} into \eqref{eq_TH32_6}, we obtain
\begin{equation*}
\|\Delta\bar\phi^n\|^2+\big{(}\frac{c_\delta}{2}-\frac{40\epsilon c_1}{\delta}\big{)}\sum\limits_{k=1}^{n-1}\tau_k\| \Delta^2\bar \phi^k\|^2
\leq \|\Delta\bar\phi^{0}\|^2+\frac{40}{\delta}(C(\epsilon,M)t_n+2(1+r_{\max})^3\epsilon\tau_1\|\Delta^2\phi^{0}\|^2).
\end{equation*}
The proof is completed.
\end{proof}

Choosing $\epsilon= \frac{c_\delta\delta}{80c_1}$ in  \eqref{eq_TH32_1}, we have the following corollary.
\begin{corollary}\label{Coroll}
  Under the condition of Theorem \ref{TH3.2_bound} for any $1\leq n\leq N$, it holds
\begin{equation*}
\|\bar{\phi}^n\|_{H^2}\leq \tilde{M}_0  \quad \text{and}\quad \|\phi^n\|_{H^2}=\|\eta^n\bar{\phi}^n\|_{H^2}\leq \gamma^0(2+\gamma^0)\tilde{M}_0:= \tilde{M},
\end{equation*}
where $\tilde{M}_0, \tilde{M}$ depend on $\Phi^0, \Omega,\delta,r_{\max}$ and $t_n$.
Noticing that $H^2(\Omega)\hookrightarrow L^\infty(\Omega)$, there exists a function $\bar{g}:\mathbb{R}\rightarrow\mathbb{R}$  such that
\begin{equation*}
  \|f(B\phi^{n-1})\|_{L^\infty},\|f'(B\phi^{n-1})\|_{L^\infty},\|f''(B\phi^{n-1})\|_{L^\infty}\leq \bar{g}(\tilde{M}),
\end{equation*}
where $\bar{g}(x)=[(1+2r_{\max})x]^3+3[(1+2r_{\max})x]^2+6[(1+2r_{\max})x]+1.$
\end{corollary}

\section{Convergence analysis}
We now consider error estimate of scheme \eqref{scheme}.
To the end, we introduce a generic positive constant $\mathcal{R}$ such that
\begin{equation}\label{Rb}
 \|\Phi\|_{L^{\infty}(0,T;H^m_{per}(\Omega))}
            +\underbrace{\|\partial_t\Phi\|_{L^\infty(0,T;H^3_{per}(\Omega))}+\|\partial_{tt}\Phi\|_{L^\infty(0,T;H^3_{per}(\Omega))}}_{\text{used in Lemma} \;\ref{Rfj}}
            +\underbrace{\|\partial_{ttt}\Phi\|_{L^\infty(0,T;H^1_{per}(\Omega))}}_{\text{used in Lemma} \; \ref{Rtj}}+1\leq \mathcal{R},
\end{equation}
where $m>3$ and
\begin{equation*}
H^m_{per}=\{u\in H^m(\Omega): u^{(k)}(0,\cdot)=u^{(k)}(L,\cdot), u^{(k)}(\cdot,0)=u^{(k)}(\cdot,L), k=0,\cdots, m-1\}.
\end{equation*}

\begin{lemma}
Assume $\bar\phi^n,\Phi(\cdot,t_n)$ are the solutions of \eqref{scheme}and \eqref{CHeq}, respectively. There exists a constant $c_\Omega$ such that
\begin{equation}\label{norm_equal}
\|\bar\phi^n-\Phi(\cdot,t_n)\|_{H^1}\leq c_\Omega\|\nabla (\bar\phi^n-\Phi(\cdot,t_n))\|.
\end{equation}
\end{lemma}
\begin{proof}
 Noting \eqref{Conser} and the initial conditions $\bar\phi^0=P_N \Phi^0$, we have
\begin{equation*}
(\bar\phi^n,1)=(\bar\phi^0,1)=(P_N\Phi^0,1)=(\Phi^0,1)=(\Phi(\cdot,t_n),1),
\end{equation*}
where the last equality follows from the mass conservation of \eqref{CHeq}.
Applying Poincar\'e-Friedrichs  inequality \cite{MR1974504}, one immediately has \eqref{norm_equal}. The proof is completed.
\end{proof}
\begin{lemma}[Discrete Gr\"{o}nwall's inequality \cite{thomee_galerkin_2006}] \label{Gronwall}
Assume that $w_n, n\geq0$ satisfy
\begin{equation*}
  w_n\leq \alpha_n+\sum\limits_{k=0}^{n-1}\beta_k w_k, \qquad\forall n\geq 0,
\end{equation*}
where $\alpha_n$ is nondecreasing and $\beta_n\geq0$. Then $w_n\leq \alpha_n\exp(\sum_{k=0}^{n-1}\beta_k).$
\end{lemma}

{We now show the following approximation results.}
\begin{lemma}[\!\!\cite{shen_spectral_2011,kreiss_stability_1979}]\label{projection_err}
For any $u\in H^m_{per}(\Omega)$ and $0\leq \lambda\leq m$, there holds
\begin{align*}
  \|P_N u-u\|_{H^\lambda_{per}}\leq C_p N^{\lambda-m}\|u\|_{H^m_{per}},
\end{align*}
where the $L^2$-orthogonal projection operator $P_N$ is defined in \eqref{OpeaP}.
\end{lemma}
Besides, the operator $P_N$ commutes with the derivation on $H^m_{per}(\Omega)$ \cite{shen_spectral_2011}, that is
\begin{equation*}
P_N\Delta^2 u=\Delta^2 P_N u, \quad \forall u\in H^4_{per}(\Omega).
\end{equation*}
Furthermore, it also holds that
\begin{equation}\label{DeltaPi}
(\Delta^2(P_N u-u), v_N)=(P_N \Delta^2 u- \Delta^2 u, v_N)=0, \quad \forall v_N\in S_N, u\in H_{per}^4(\Omega).
\end{equation}
We further introduce the discrete complementary convolution (DCC) kernel $p_{n-j}^{(n)}$ \cite{zhang_sharp_2020} defined by
\begin{equation}
\sum_{j=k}^n p_{n-j}^{(n)} b_{j-k}^{(j)} \equiv1, \quad  \forall  1\leq k \leq n,  \; 1 \leq n \leq N.\label{defin_p}
\end{equation}
The DOC and DCC  kernels have the following  relationship (c.f. \cite[Proposition 2.1]{zhang_sharp_2020})
\begin{equation*}
p^{(n)}_{n-j}=\sum_{l=j}^n \theta^{(l)}_{l-j},\quad \theta^{(n)}_{n-j}= p^{(n)}_{n-j}-p^{(n-1)}_{n-1-j},\quad \forall 1\leq j\leq n,
\end{equation*}
where $p^{(n)}_{-1}:=0 \; (\forall n\geq 0)$ is defined. The DCC kernel is valuable in the global error estimate \cite{zhang_sharp_2020}.
\begin{proposition}[\!\!\cite{zhang_sharp_2020}]\label{Proposition_Pestimate}
The DCC kernels $p^{(n)}_{n-k}$ defined in \eqref{defin_p} have the following properties
\begin{equation*}
\sum_{j=1}^n p^{(n)}_{n-j} = t_n, \quad p^{(n)}_{n-j}\leq 2\tau \quad \text{for any} \quad 1\leq j\leq n.
\end{equation*}
\end{proposition}

The following lemma focuses on the local and global consistency errors.
\begin{lemma}\label{Rfj}
Set $R_f^j= f(B\Phi(\cdot,t_{n-1}))-f(\Phi(\cdot,t_n))$. It holds that
\begin{align}
&\|\nabla\Delta R_f^1\|\leq 15\bar{g}(\mathcal{R})\mathcal{R}^3\tau, \,
 \|\nabla\Delta R_f^k\|\leq 15\bar{g}(\mathcal{R})\mathcal{R}^3(1+5r_{\max})t_n\tau^2, k>1,\label{Tru1}\\
 &\|\nabla R_f^1\|\leq 2\bar{g}(\mathcal{R})\mathcal{R}^2\tau,\,\,\quad
 \|\nabla R_f^k\|\leq 2\bar{g}(\mathcal{R})\mathcal{R}^2(1+5r_{\max})t_n\tau^2, k>1. \label{Tru2}
\end{align}
Furthermore,
\begin{equation}\label{two_ord}
\sum\limits_{k=1}^n\|\sum\limits_{j=1}^k\theta^{(k)}_{k-j}\nabla\Delta R_f^j\|\leq 15\bar{g}(\mathcal{R})\mathcal{R}^2[t_n(1+5r_{\max})+2]\tau^2.
\end{equation}
\end{lemma}
\begin{proof}
By the mean value theorems, we have
$$R_f^j=f'((1-\lambda)B\Phi(\cdot,t_{j-1})+\lambda\Phi(\cdot,t_j)) (B\Phi(\cdot,t_{j-1})-\Phi(\cdot,t_j)):=f'(\Theta)R, $$
where $0\leq \lambda\leq 1$, and $\Theta, R\in H^m_{per}(\Omega)$. A direct calculation gives
\begin{align*}
\nabla\Delta (f'(\Theta)R)=&\underbrace{f''(\Theta)\nabla\Delta \Theta R+f'''(\Theta)\Delta \Theta\nabla \Theta R+f'''(\Theta)\nabla|\nabla \Theta|^2R+0}_{\nabla\Delta f'(\Theta) R}\\
&+\underbrace{f''(\Theta)\Delta \Theta \nabla R +f'''(\Theta)|\nabla \Theta|^2\nabla R}_{\Delta f'(\Theta)\nabla R}+\underbrace{f''(\Theta)\nabla \Theta\Delta R}_{\nabla f'(\Theta)\Delta R}+f'(\Theta)\nabla\Delta R\\
&+\underbrace{2f''(\Theta)\nabla(\nabla \Theta\cdot\nabla R)+2f'''(\Theta)\nabla \Theta (\nabla \Theta\cdot\nabla R)}_{2\nabla(\nabla f'(\Theta)\cdot\nabla R)},
\end{align*}
where $f''''(\Theta)=0$ has been used.
In view of
\begin{align*}
\|\nabla(\nabla \Theta\cdot\nabla R)\|&=\|\nabla(\Theta_xR_x+\Theta_yR_y)\|\leq\|\Theta_xR_x+\Theta_yR_y\|_{H^1}\\
&\leq \|\Theta_xR_x\|_{H^1}+\|\Theta_yR_y\|_{H^1}\leq \|\Theta_x\|_{W^{1,4}}\|R_x\|_{W^{1,4}}+\|\Theta_y\|_{W^{1,4}}\|R_y\|_{W^{1,4}}\\
&\leq 2\|\Theta\|_{W^{2,4}}\|R\|_{W^{2,4}},
\end{align*}
it follows from H\"{o}lder inequality and Sobolev embedding theorems that
\begin{align}
\|\nabla\Delta (f'(\Theta)R)\|\leq &\|f''(\Theta)\|_{L^\infty}\|\Theta\|_{H^3}\| R\|_{L^\infty}+\|f'''(\Theta)\|_{L^\infty}\|\Theta\|_{W^{2,3}}\|\Theta \|_{W^{1,6}}\|R\|_{L^\infty}\nonumber\\
&+\|f'''(\Theta)\|_{L^\infty}(2\|\Theta\|_{W^{2,4}}^2)\|R\|_{L^\infty}\nonumber\\
&+\|f''(\Theta)\|_{L^\infty}\|\Theta\|_{W^{2,3}} \|R\|_{W^{1,6}} +\|f'''(\Theta)\|_{L^\infty}\|\Theta\|_{W^{1,4}}^2\|R\|_{W^{1,\infty}}\nonumber\\
&+\|f''(\Theta)\|_{L^\infty}\| \Theta\|_{W^{1,6}}\|R\|_{W^{2,3}}+\|f'(\Theta)\|_{L^\infty}\|R\|_{H^3}\nonumber\\
&+4\|f''(\Theta)\|_{L^\infty}\|\Theta\|_{W^{2,4}}\|R\|_{W^{2,4}}
+4\|f'''(\Theta)\|_{L^\infty}\|\Theta\|_{W^{1,6}}^2\| R\|_{W^{1,6}}\nonumber\\
\leq& M_f^j\|B\Phi(\cdot,t_{j-1})-\Phi(\cdot,t_j)\|_{H^3},
\label{L43_1}
\end{align}
where $M_f^j=15C_\Omega\max\{\|f'(\Theta)\|_{L^\infty},\|f''(\Theta)\|_{L^\infty},\|f'''(\Theta)\|_{L^\infty}\}\max\{\|\Theta\|_{H^3},\|\Theta\|_{H^3}^2\}\leq 15C_\Omega\bar{g}(\mathcal{R})\mathcal{R}^2$.
By Taylor expansion, for $j\geq2$, we have
\begin{equation}\label{taylor}
  B\Phi(\cdot,t_{j-1})-\Phi(\cdot,t_j)=(1+r_n)\int_{t_{n-1}}^{t_n}(s-t_{n-1})\partial_{tt}\Phi(s)ds
-r_n\int_{t_{n-2}}^{t_n}(s-t_{n-2})\partial_{tt}\Phi(s)ds.
\end{equation}
It is a consequence of \eqref{taylor} that
\begin{equation}
 \|B\Phi(\cdot,t_{j-1})-\Phi(\cdot,t_j)\|_{H^3}\leq (1+5r_{\max})\tau^2\|\partial_{tt}\Phi\|_{L^\infty(0,T;H^3(\Omega))},\quad j\geq 2.
\label{L43_2}
\end{equation}
Similarly, for $j=1$, we have
\begin{equation}
  \|B\Phi^0-\Phi(\cdot,t_1)\|_{H^3}\leq \tau\|\partial_{t}\Phi\|_{L^\infty(0,T;H^3(\Omega))}.
\label{L43_3}
\end{equation}
Combining \eqref{L43_1}-\eqref{L43_3}, we can deduce \eqref{Tru1}.
A similar argument to \eqref{Tru1} will show \eqref{Tru2}.

The result of \eqref{two_ord} can be proved by the Proposition \ref{Proposition_Pestimate} as follows.
\begin{align*}
\sum\limits_{k=1}^n\|\sum\limits_{j=1}^k\theta^{(k)}_{k-j}\nabla\Delta R_f^j\|
&\leq \sum\limits_{k=1}^n\sum\limits_{j=1}^k\theta^{(k)}_{k-j}\|\nabla\Delta R_f^j\|
\leq \sum\limits_{j=1}^n\|\nabla\Delta R_f^j\|\sum\limits_{k=j}^n\theta^{(k)}_{k-j}=\sum\limits_{j=1}^n\|\nabla\Delta R_f^j\|p^{(n)}_{n-j}\\
&\leq \max_{2\leq j\leq n}\|\nabla\Delta R_f^j\| \sum\limits_{j=2}^np^{(n)}_{n-j}+\|\nabla\Delta R_f^1\|p^{(n)}_{n-1} \leq 15\bar{g}(\mathcal{R})\mathcal{R}^2[t_n(1+5r_{\max})+2]\tau^2.
\end{align*}
The proof is completed.
\end{proof}
\begin{lemma}[\!\!\cite{zhang_sharp_2020}] \label{Rtj}
Denote $R_t^j=\mathcal{D}_2 \Phi(\cdot,t_j)- \partial_t\Phi (t_j).$ It holds
\begin{equation*}
  \|R_t^1\|_{H^1}\leq \tau\|\partial_{tt}\Phi\|_{L^\infty(0,\tau;H^1(\Omega))},\quad
  \|R_t^j\|_{H^1}\leq \frac{1+5r_{\max}}{2} \tau^2\|\partial_{ttt}\Phi\|_{L^\infty(0,T;H^1(\Omega))}, \; j>1 .
\end{equation*}
Moreover, the global truncation error may be estimated by
\begin{equation*}
\sum\limits_{k=1}^n\|\sum\limits_{j=1}^k\theta^{(k)}_{k-j}R_t^j\|_{H^1}\leq (2+\frac{1+5r_{\max}}{2}t_n)\tau^2\mathcal{R}.
\end{equation*}
\end{lemma}
The proof of Lemma \ref{Rtj} is similar to that of Lemma \ref{Rfj} and we omitted it here. We denote
\begin{align*}
&s^n=\gamma^n-\Gamma(t_n);                   &&\bar e_N^n=\bar \phi^n- P_N \Phi(\cdot,t_n);       && e_N^n= \phi^n-P_N \Phi(\cdot,t_n); \nonumber\\
&e_P^n= P_N \Phi(\cdot,t_n)-\Phi(\cdot,t_n); &&\bar e^n  =\bar \phi^n- \Phi(\cdot,t_n);    &&e^n=\phi^n-\Phi(\cdot,t_n).
\end{align*}

{To simplify the notations, we further denote the maximums of sequence as}
\begin{equation}\label{Dnote}
  \|\nabla\bar e_N^{n_0}\|=\max\limits_{0\leq k \leq n}\|\nabla\bar e_N^k\|; \quad \rho^n=\max\limits_{0\leq k \leq n}\{\|e_P^k\|_{H^1} +|1-\eta^k|\|\bar\phi^k\|_{H^1}\}/c_\Omega.
\end{equation}
\begin{lemma}\label{LemEFk}
Set $E_f^k=f(B\phi^{k-1})-f(B\Phi(\cdot,t_{k-1}))$. It holds
\begin{align}
&\|\nabla E_f^k\|\leq \bar{C}_{1}(\|\nabla\bar e_N^{n_0}\|+\rho^{n-1}),\, 1\leq k\leq n, \label{L46_1}\\
&\sum\limits_{k=1}^n\tau_k\|\nabla E_f^k\|^2
\leq  \bar{C}_{2}(\|\nabla\bar e_N^{n_0}\|+\rho^{n-1})(t_n\rho^{n-1}+\sum\limits_{k=1}^{n}\tau_k\|\nabla \bar e_N^{k-1}\|),
\label{L46_2}
\end{align}
where $\bar{C}_{1}=\bar{g}(\tilde{M})\bar{g}(\tilde{M}+\mathcal{R})(1+ 2r_{\max})c_\Omega$, $\bar{C}_{2}= 2(1+ r_{\max})^3[\bar{g}(\tilde{M})\bar{g}(\tilde{M}+\mathcal{R})c_\Omega]^2$ with $c_\Omega$ given in \eqref{norm_equal}.
\end{lemma}
\begin{proof}
In view of $\nabla g(u)=g'(u)\nabla u$, we have
\begin{align*}
\nabla E_f^k&=f'(B\phi^{k-1})\nabla B\phi^{k-1}-f'(B\Phi(\cdot,t_{k-1}))\nabla B\Phi(\cdot,t_{k-1})\\
&=f'(B\phi^{k-1})\nabla Be^{k-1}+(f'(B\phi^{k-1})-f'(B\Phi(\cdot,t_{k-1})))\nabla B\Phi(\cdot,t_{k-1})\\
&=f'(B\phi^{k-1})\nabla Be^{k-1}+f''(\Theta)Be^{k-1}\nabla B\Phi(\cdot,t_{k-1}),
\end{align*}
where $\Theta=(1-\lambda)B\phi^{j-1}+\lambda B\Phi(\cdot,t_{j-1}), 0\leq \lambda\leq 1$.
From Corollary \ref{Coroll}, we get
\begin{align}
\|\nabla E_f^k\|&\leq \bar{g}(\tilde{M})\|\nabla B e^{k-1}\|+\bar{g}(\tilde{M})\|Be^{k-1}\|_{L^6}\|\nabla B\Phi(\cdot,t_{k-1})\|_{L^3}\nonumber\\
&\leq C_\Omega\bar{g}(\tilde{M})\bar{g}(\tilde{M}+\mathcal{R})\|B e^{k-1}\|_{H^1}:=C_{\tilde{M},\mathcal{R}}\|B e^{k-1}\|_{H^1},
\label{eq_L46_1}
\end{align}
where $\|\nabla B\Phi(\cdot,t_{k-1})\|_{L^3}\leq C_\Omega\|B\Phi(\cdot,t_{k-1})\|_{H^2}\leq \mathcal{R}$ and $\|Be^{k-1}\|_{L^6}\leq C_\Omega\|Be^{k-1}\|_{H^1}$ have been used.
Similar to \eqref{eq_TH32_8}, we can obtain
\begin{equation}
\sum_{k=1}^n\tau_k\|B e^{k-1}\|_{H^1}^2\leq  2(1+ r_{\max})^3 \sum_{k=1}^{n}\tau_k\|e^{k-1}\|_{H^1}^2.
\end{equation}
It follows from the definitions in \eqref{Dnote} and \eqref{scheme4} that
\begin{align}
\|e^{k-1}\|_{H^1}&\leq  \|\bar e_N^{k-1}\|_{H^1} + \| e_P^{k-1}\|_{H^1} +|1-\eta^{k-1}|\|\bar\phi^{k-1}\|_{H^1}\nonumber\\
&\leq c_\Omega(\|\nabla \bar e_N^{k-1}\|+\rho^{k-1}) \leq c_\Omega(\|\nabla\bar e_N^{n_0}\|+\rho^{k-1}),
\label{eq_L46_3}
\end{align}
where the second inequality used \eqref{norm_equal}.
Thus, we have
\begin{align}
\sum\limits_{k=1}^n\tau_k\|e^{k-1}\|_{H^1}^2&\leq c_\Omega^2(\|\nabla\bar e_N^{n_0}\|+\rho^{n-1})\sum\limits_{k=1}^{n}\tau_k(\|\nabla \bar e_N^{k-1}\|+\rho^{k-1})\nonumber\\
&\leq c_\Omega^2(\|\nabla\bar e_N^{n_0}\|+\rho^{n-1})(t_n\rho^{n-1}+\sum\limits_{k=1}^{n}\tau_k\|\nabla \bar e_N^{k-1}\|).
\label{eq_L46_4}
\end{align}
It follows from \eqref{eq_L46_1}-\eqref{eq_L46_4} that
\begin{align*}
\sum\limits_{k=1}^n\tau_k\|\nabla E_f^k\|^2
& \leq 2(1+ r_{\max})^3C_{\tilde{M},\mathcal{R}}^2 \sum_{k=1}^{n}\tau_k\|e^{k-1}\|_{H^1}^2 \nonumber\\
&\leq 2(1+ r_{\max})^3C_{\tilde{M},\mathcal{R}}^2c_\Omega^2(\|\nabla\bar e_N^{n_0}\|+\rho^{n-1})(t_n\rho^{n-1}+\sum\limits_{k=1}^{n}\tau_k\|\nabla \bar e_N^{k-1}\|).
\end{align*}
which implies \eqref{L46_2}.

On the other hand, noticing that $\|B e^{k-1}\|_{H^1}\leq (1+2r_{\max})\max\limits_{1\leq k\leq n}\|e^{k-1}\|_{H^1}$ and \eqref{eq_L46_1}, \eqref{eq_L46_3}, we can deduce \eqref{L46_1}.
The proof is completed.
\end{proof}

Now we give the following error estimate in $H^1$-norm.
\begin{theorem}\label{Th5}
Assume $\Phi(\cdot,t_n)\in L^\infty(0,T;H^m_{per}(\Omega)), (m>3)$ and $\phi^n, \bar\phi^n\in S_N$ are the solutions of \eqref{CHeq} and \eqref{scheme}, respectively.
If
\begin{equation}\label{Limit_tau}
\tau < \min\{\frac{1}{(1+\bar{g}(\mathcal{R}))^2\mathcal{R}^2+1}, \frac{1}{2 (c_0^2+1)}\},\quad  N^{3-m}\leq \frac{1}{2 (c_0^2+1)},
\end{equation}then
\begin{align}
&\|\nabla e^n\|, \|\nabla \bar e^n\| \leq C(\tau^2+N^{3-m}),\label{re1}\\
&\sum\limits_{k=1}^n\tau_k\|\nabla\Delta e^k\|^2, \sum\limits_{k=1}^n\tau_k\|\nabla\Delta \bar e^k\|^2 \leq C(\tau^4+N^{2(3-m)}),\label{re2}\\
&s^n\leq C(\tau+N^{3-m}),\label{re3}
\end{align}
where constant $C>0$ is independent of $\tau$ and $N$, and $c_0$ will be given in the proof.
\end{theorem}
\begin{proof}
We now prove the theorem by mathematical induction. It is easy to check that $\xi^0=1$ and it holds $|1-\xi^0|=0\leq c_0(\tau+N^{3-m})$.
Suppose that it holds for any $0\leq q\leq n-1$ that
\begin{equation}\label{suppo}
|1-\xi^q|\leq c_0(\tau+N^{3-m}),
\end{equation}
then the mathematical induction is completed if we derive
\begin{equation}\label{conclution}
|1-\xi^n|\leq c_0(\tau+N^{3-m}),\quad \forall n\leq K.
\end{equation}
where $c_0$ will be determined below.

For clarity, the proof of \eqref{conclution} is  divided into the following 5 steps.

\textbf{ Step 1: The boundeness of} $\rho^n$.
By the definition of $\eta$ in \eqref{scheme4} and \eqref{suppo},  we have
\begin{equation}
|1-\eta^q|=|1-\xi^q|^2\leq c_0^2(\tau+N^{3-m})^2.
\label{tri2}
\end{equation}
Combining \eqref{Limit_tau} and \eqref{suppo}, one has $|1-\xi^q|\leq \frac12$ and
$|\eta^q|\geq \frac34$. It follows from \eqref{scheme4} and Theorem \ref{TH3.1} that
$$
\|\bar\phi^q\|_{H^1}= \frac{\|\phi^q\|_{H^1}}{|\eta^q|}\leq \frac{4M}{3}.
$$
Consequently, with the definition of $\rho^n$ in \eqref{Dnote},  we have
\begin{equation}\label{term2.1}
\rho^{n-1}\leq \frac{C_p\mathcal{R}}{c_\Omega}N^{1-m}+\frac{4M}{3c_\Omega}c_0^2(\tau+N^{3-m})^2
\leq (\frac{C_p\mathcal{R}}{c_\Omega}+\frac{4M}{3c_\Omega})(N^{1-m}+c_0^2(\tau+N^{3-m})^2),
\end{equation}
where Lemma \ref{projection_err} has been used.

\textbf{Step 2: Estimate for }$\|\nabla\bar e_N^n\|$.
From \eqref{scheme1}, we have the error equation
\begin{equation}\label{error_equation}
(\D 2 \bar e^j,v_N)+(\Delta\bar e^j, \Delta v_N)-(\Delta E_f^j,v_N)=(\Delta R_{f}^j,v_N)+(R_t^j,v_N),
\end{equation}
where $E_f^j, R_{f}^j, R_t^j$ are given by
\begin{equation*}
E_f^j=f(B\phi^{j-1})- f(B\Phi(\cdot,t_{j-1})), \
R_f^j=f(B\Phi(\cdot,t_{j-1}))-f(\Phi(\cdot,t_j)), \
R_t^j=\partial_t\Phi(\cdot,t_j)-\D 2 \Phi(\cdot,t_j).
\end{equation*}
By the definition of  $P_N$  and the properties in \eqref{DeltaPi}, one has
$$ (\D 2 e_P^j, v_N)=0 \quad \text{and} \quad  (\Delta^2 e_P^j,  v_N)=0, \quad \forall v_N\in S_N.$$
 Thus \eqref{error_equation} can be rewritten as
\begin{equation}\label{error_equation2}
(\D 2 \bar e_N^j,v_N)+(\Delta^2\bar e_N^j, v_N)-(\Delta E_f^j,v_N)=(\Delta R_{f}^j,v_N)+(R_t^j,v_N).
\end{equation}
Multiplying \eqref{error_equation2} by $\theta^{(k)}_{k-j}$ and summing $j$ from
1 to $k$, we have
\begin{equation*}
(\nabla_\tau\bar e_N^{k},v_N)
+\sum\limits_{j=1}^k\theta^{(k)}_{k-j}(\Delta^2\bar e_N^j,v_N)
-\sum\limits_{j=1}^k\theta^{(k)}_{k-j}(\Delta E_f^{j},v_N)
=\sum\limits_{j=1}^k\theta^{(k)}_{k-j}(\Delta R_{f}^j+R_t^j,v_N),
\end{equation*}
where \eqref{Theta_D2} has been used.
Choosing $v_N=-2\Delta\bar  e_N^{k}$ and summing $k$ from
1 to $n$, we have
\begin{align}\label{e462}
\|\nabla\bar e_N^n\|^2 &-\|\nabla\bar e_N^0\|^2+2\sum\limits_{k=1}^n(\nabla_\tau\nabla\bar e_N^k)^2
+2\sum\limits_{k=1}^n\sum\limits_{j=1}^k\theta^{(k)}_{k-j}(\nabla \Delta\bar e_N^j,\nabla\Delta\bar e_N^k)\nonumber\\
&\qquad =2\sum\limits_{k=1}^n\sum\limits_{j=1}^k\theta^{(k)}_{k-j}(\nabla E_f^{j},\nabla\Delta\bar e_N^k)
+2\sum\limits_{k=1}^n\sum\limits_{j=1}^k\theta^{(k)}_{k-j}(\nabla\Delta R_{f}^j+\nabla R_t^j,\nabla\bar e_N^k).
\end{align}
Exchanging the summation order and applying Young's inequality, we have
\begin{align}
2\sum\limits_{k=1}^n\sum\limits_{j=1}^k\theta^{(k)}_{k-j}(\nabla E_f^{j},\nabla \Delta\bar e_N^k)
&=2\sum\limits_{j=1}^n(\nabla E_f^{j},\sum\limits_{k=j}^n\theta^{(k)}_{k-j}\nabla \Delta\bar e_N^k)\nonumber\\
&\leq 40\delta^{-1}\sum\limits_{j=1}^n\tau_j\|\nabla E_f^j\|^2
+\delta\sum\limits_{j=1}^n\frac{\|\sum_{k=j}^n\theta^{(k)}_{k-j}\nabla \Delta\bar e_N^k\|^2}{40\tau_j}.
\label{eq453}
\end{align}
From Lemma \ref{Lem_inequal}, we have
\begin{equation}
2\sum\limits_{k=1}^n\sum\limits_{j=1}^k\theta^{(k)}_{k-j}(\nabla \Delta\bar e_N^j,\nabla\Delta\bar e_N^k)
\geq \delta\sum\limits_{j=1}^n\frac{\|\sum_{k=j}^n\theta^{(k)}_{k-j}\nabla \Delta\bar e_N^k\|^2}{20\tau_j}.
\label{eq454}
\end{equation}
Inserting \eqref{eq453}-\eqref{eq454} into \eqref{e462} and
removing $2\sum\limits_{k=1}^n(\nabla_\tau\nabla\bar e_N^k)^2$ yields
\begin{align}
&\|\nabla\bar e_N^n\|^2 +\delta\sum\limits_{j=1}^n\frac{\|\sum_{k=j}^n\theta^{(k)}_{k-j}\nabla \Delta\bar e_N^k\|^2}{40\tau_j} \nonumber\\
\leq & \|\nabla\bar e_N^0\|^2+40\delta^{-1}\sum\limits_{j=1}^n\tau_j\|\nabla E_f^j\|^2 +2\sum\limits_{k=1}^n\|\sum\limits_{j=1}^k\theta^{(k)}_{k-j}(\nabla \Delta R_{f}^j+\nabla R_{t}^j)\|\|\nabla\bar e_N^k\|.
\label{err_1}
\end{align}

We now estimate \eqref{err_1} item by item.
It follows from  Lemmas \ref{Rfj} and \ref{Rtj} that
\begin{align}\label{trunerr}
&2\sum\limits_{k=1}^n\|\sum\limits_{j=1}^k\theta^{(k)}_{k-j}\nabla \Delta R_{f}^j\|+
2\sum\limits_{k=1}^n\|\sum\limits_{j=1}^k\theta^{(k)}_{k-j}\nabla R_t^j\|\nonumber\\
 \leq  & 30\bar{g}(\mathcal{R})\mathcal{R}^2[t_n(1+5r_{\max})+2]\tau^2+(4+(1+5r_{\max})t_n)\mathcal{R}\tau^2:=\bar{C}_{3}\tau^2.
\end{align}
  If $\|\nabla\bar e_N^{n_0}\|<\rho^{n-1}$,  one can directly obtain the following estimate by  \eqref{term2.1}
\begin{equation}\label{eBrho}
\|\nabla\bar e_N^{n_0}\|\leq (\frac{C_p\mathcal{R}}{c_\Omega}+\frac{4M}{3c_\Omega})(N^{1-m}+c_0^2(\tau+N^{3-m})^2).
\end{equation}
Now we only need to consider the situation of $\|\nabla\bar e_N^{n_0}\|\geq\rho^{n-1}$.
It follows from Lemma \ref{LemEFk} that
\begin{align}
40\delta^{-1}\sum\limits_{k=1}^n\tau_k\|\nabla E_f^k\|^2
&\leq 40\delta^{-1}\bar{C}_{2}(\|\nabla\bar e_N^{n_0}\|+\rho^{n-1})(t_n\rho^{n-1}+\sum\limits_{k=1}^{n}\tau_k\|\nabla \bar e_N^{k-1}\|)\label{est_E2}\\
&\leq 80\delta^{-1}\bar{C}_{2}\|\nabla\bar e_N^{n_0}\|(t_n\rho^{n-1}+\sum\limits_{k=1}^{n}\tau_k\|\nabla \bar e_N^{k-1}\|).
\label{est_E}
\end{align}
Here $\|\nabla\bar e_N^{n_0}\|$ is used instead of \eqref{term2.1} to estimate $\rho^{n-1}$.
This idea of inequality zoom $\rho^{n-1}\leq \|\nabla\bar e_N^{n_0}\|$ will benefit our later estimation.

Inserting \eqref{trunerr} and \eqref{est_E} into \eqref{err_1} and using the definition of $\|\nabla\bar e_N^{n_0}\|$ , we have
\begin{align}
&\|\nabla\bar e_N^n\|\|\nabla\bar e_N^{n_0}\|\leq\|\nabla\bar e_N^{n_0}\|^2\nonumber\\
\leq &\|\nabla\bar e_N^0\|^2+\bar{C}_{3}\tau^2\|\nabla\bar e_N^{n_0}\|+80\delta^{-1}\bar{C}_{2}\|\nabla\bar e_N^{n_0}\|
(t_{n_0}\rho^{n_0-1}+\sum\limits_{k=1}^{n_0}\tau_k\|\nabla \bar e_N^{k-1}\|)\nonumber\\
\leq &\|\nabla\bar e_N^0\|\|\nabla\bar e_N^{n_0}\|+\bar{C}_{3}\tau^2\|\nabla\bar e_N^{n_0}\|+80\delta^{-1}\bar{C}_{2}\|\nabla\bar e_N^{n_0}\|
(t_n\rho^{n-1}+\sum\limits_{k=1}^{n}\tau_k\|\nabla \bar e_N^{k-1}\|).\label{err}
\end{align}
Eliminating $\|\nabla\bar e_N^{n_0}\|$ from \eqref{err} gives
\begin{equation*}
\|\nabla\bar e_N^n\|
\leq  \|\nabla\bar e_N^0\|+\bar{C}_{3}\tau^2
+80\delta^{-1}\bar{C}_{2}(t_n\rho^{n-1}+\sum\limits_{k=1}^{n}\tau_k\|\nabla \bar e_N^{k-1}\|).
\end{equation*}
By Lemma \ref{Gronwall},  we have
\begin{equation}
\|\nabla \bar e_N^n\|\leq \exp(80\delta^{-1}\bar{C}_{2} t_n)
\left(\|\nabla\bar e_N^0\|+\bar{C}_{3}\tau^2+80\delta^{-1}\bar{C}_{2}t_n\rho^{n-1}\right).
\end{equation}
Noticing that  $\phi^0=P_N \Phi^0$ and \eqref{term2.1}, we get
\begin{equation}
\|\nabla \bar e_N^n\|\leq C_0\left(\tau^2+N^{1-m}+c_0^2(\tau+N^{3-m})^2\right),\label{e4660}
\end{equation}
where $$C_0=\exp(80\delta^{-1}\bar{C}_{2} T)\left[C_p\mathcal{R}+\bar{C}_{3}+80\delta^{-1}\bar{C}_{2}T(\frac{C_p\mathcal{R}}{c_\Omega}+\frac{4M}{3c_\Omega})\right].$$
Combining \eqref{e4660} and \eqref{eBrho}, we have the following estimation whether $\|\nabla\bar e_N^{n_0}\|> \rho^{n-1}$ or not
\begin{equation}
 \|\nabla \bar e_N^n\|\leq \max\{\|\nabla \bar e_N^n\|, \rho^{n-1}\}\leq C_1\left(\tau^2+N^{1-m}+c_0^2(\tau+N^{3-m})^2\right),\label{e466}
\end{equation}
where $C_1=C_0+(\frac{C_p\mathcal{R}}{c_\Omega}+\frac{4M}{3c_\Omega})$.
According to Lemma \ref{Lem_inequal} and \eqref{err_1}-\eqref{est_E2}, we have
\begin{align}
\frac{c_\delta}{2}\sum_{k=1}^n\tau_k\|\nabla \Delta \bar e_N^n\|^2
&\leq \|\nabla \bar e_N^{n_0}\|^2+40\delta^{-1}\bar{C}_{2}t_n(\|\nabla \bar e_N^{n_0}\|+\rho^{n-1})^2+\bar{C}_3\tau^2\|\nabla \bar e_N^{n_0}\|
.\label{e467}
\end{align}

\textbf{Step 3: Give an estimate of} $\|\nabla \bar{ e}^n\|$ \textbf{without introducing} $c_0$.
It follows from \eqref{Limit_tau} that
\begin{equation}
  \tau^2+N^{1-m}+c_0^2(\tau+N^{3-m})^2\leq 2(\tau+N^{3-m}).\label{simplify}
\end{equation}
Considering \eqref{e466}-\eqref{e467}, Lemma \ref{projection_err} and the triangle inequality yields
\begin{align}
\|\nabla \bar e^n\|&\leq\|\nabla\bar e_N^n\|+\| e_P^n\|_1\nonumber\\
&\leq C_1(\tau^2+N^{1-m}+c_0^2(\tau+N^{3-m})^2)+C_p\mathcal{R}N^{1-m}\label{result12}\\
&\leq C_2(\tau+N^{3-m}),\label{result1}
\end{align}
and
\begin{align}
\sum_{k=1}^n\tau_k\|\nabla \Delta \bar e^k\|^2&\leq2\sum_{k=1}^n\tau_k(\|\nabla\Delta \bar e_N^k\|^2+\|\nabla\Delta e_P^k\|^2)\nonumber\\
&\leq\frac{4}{c_\delta}\Big[\|\nabla \bar e_N^{n_0}\|^2+40\delta^{-1}\bar{C}_{2}t_n(\|\nabla \bar e_N^{n_0}\|+\rho^{n-1})^2+\bar{C}_3\tau^2\|\nabla \bar e_N^{n_0}\|\Big]+2t_n\|\nabla\Delta e_P^{n_0}\|^2\nonumber\\
&\leq \frac{4}{c_\delta}(1+160\delta^{-1}\bar{C}_{2}t_n)C_1^2\left(\tau^2+N^{1-m}+c_0^2(\tau+N^{3-m})^2\right)^2\nonumber\\
&\quad +\frac{4}{c_\delta}\bar{C}_3\tau^2C_1\left(\tau^2+N^{1-m}+c_0^2(\tau+N^{3-m})^2\right)+2t_n(C_p\mathcal{R}N^{3-m})^2\label{result22}\\
&\leq C_2(\tau^2+N^{2(3-m)}),\label{result2}
\end{align}
where the last steps of \eqref{result1} and \eqref{result2} follow from \eqref{simplify}. Here $C_2$ is independent of $c_0$ and defined as
$$C_2=(2C_1+C_p\mathcal{R})+[\frac{32}{c_\delta}(1+160\delta^{-1}\bar{C}_{2}T)C_1^2+\frac{8}{c_\delta}\bar{C}_3\tau^2C_1+2T(C_p\mathcal{R})^2].$$
According to \eqref{L46_1}, \eqref{Limit_tau}, \eqref{term2.1} and \eqref{result1}, we have
\begin{align*}
\|\nabla E_f^k\|
&\leq \bar{C}_{1}[C_2(\tau+N^{3-m})+C_p\mathcal{R}c_\Omega ^{-1}N^{1-m}+\frac{4M}{3c_\Omega }(\tau+N^{3-m})].
\end{align*}
Combining with \eqref{Tru2} in Lemma \ref{Rfj}, we have
\begin{equation}\label{Ej}
\|\nabla E_f^k\|+\|\nabla R_f^k\|\leq C_3(\tau+N^{3-m}),
\end{equation}
where $C_3=\bar{C}_{1}[C_2+C_p\mathcal{R}c_\Omega^{-1} +\frac{4M}{3c_\Omega }]+2\bar{g}(\mathcal{R})\mathcal{R}^2[(1+5r_{\max})T^2+1].$
Besides,
\begin{align}\label{E_E}
E(\bar\phi^k)-E(\Phi(\cdot,t_k))&=\frac{1}{2}\int_\Omega |\nabla\bar\phi^k|^2-|\nabla\Phi(\cdot,t_k)|^2 {\rm{d}x} +\frac{1}{4\varepsilon^2}\int [(\bar\phi^k)^2-1]^2-[(\Phi(\cdot,t_k)^2-1]^2){\rm{d}x}\nonumber\\
&=\frac{1}{2}\big(\nabla(\bar\phi^k+\Phi(\cdot,t_k)),\nabla\bar e^k \big)
+\frac{1}{4\varepsilon^2}\big( (\bar\phi^k)^2+\Phi(\cdot,t_k)^2-2,(\bar\phi^k+\Phi(\cdot,t_k))\bar e^k\big)\nonumber\\
&\leq \bar{C}_4\|\nabla\bar e^k\|\leq \bar{C}_4C_2(\tau+N^{3-m}):= C_4(\tau+N^{3-m}),
\end{align}
where $\bar{C}_4=\frac{1}{2}(M+\mathcal{R})+\frac{1}{4\varepsilon^2}(C_\Omega^2\tilde{M}^2+\mathcal{R}^2+2)(M+\mathcal{R})$ and
the first inequality follows from $\|\bar \phi^n\|_{L^\infty}\leq C_\Omega\|\bar\phi^n\|_{H^2}\leq C_\Omega\tilde{M}$ and $\|\bar \phi^n\|_{H^1}\leq M$.

\textbf{Step 4: Estimate of }$|s^n|$.
Noticing that
\begin{align}
&\partial_t\Gamma(t_n)+\|\nabla\mu(t_n)\|^2=0, \label{4.1} \\
&\frac{\gamma^n-\gamma^{n-1}}{\tau_n}+\frac{\gamma^n}{E(\bar\phi^n)+1}\|-\nabla\Delta\bar\phi^n+\nabla f(B\phi^{n-1})\|^2=0,\label{4.2}
\end{align}
 we can get the following error equation
\begin{align*}
&\frac{s^k-s^{k-1}}{\tau_k}+\frac{\gamma^k}{E(\bar\phi^k)+1}\|-\nabla\Delta\bar\phi^k+\nabla f(B\phi^{k-1})\|^2-\|\nabla\mu(t_k)\|^2=R_\gamma^k,
\end{align*}
where $R_\gamma^k=\partial_t\Gamma(t_k)-\frac{\Gamma(t_k)-\Gamma(t_{k-1})}{\tau_k}.$
Multiplying $\tau_k$ and summing $k$ from 1 to n, we have
\begin{equation}\label{err_gamma}
s^n=s^0-\sum_{k=1}^n\tau_k\|\nabla \mu(t_k)\|^2Q_1^k-\sum_{k=1}^n\tau_k\frac{\gamma^k}{E(\bar\phi^k)+1}Q_2^k+\sum_{k=1}^n\tau_{k}R_\gamma^k,
\end{equation}
where
$Q_1^k=\frac{\gamma^k}{E(\bar\phi^k)+1}-1;\quad Q_2^k=\|-\nabla\Delta\bar\phi^n+\nabla f(B\phi^{k-1})\|^2-\|\nabla\mu(t_k)\|^2.$

We now estimate $Q_1, Q_2$ one by one. In view of $\Gamma(t_k) =E(\Phi(\cdot, t_k))+1$ and \eqref{E_E}, we have
\begin{align}\label{Q_1}
Q_1^k&=\frac{s^n}{E(\bar\phi^k)+1}+\frac{\Gamma(t_k)}{E(\bar\phi^k)+1}-1=\frac{s^n}{E(\bar\phi^k)+1}+\frac{E(\Phi(\cdot,t_k))-E(\bar\phi^k)}{E(\bar\phi^k)+1} \leq s^n+C_4(\tau+N^{3-m}),
\end{align}
In view of
$
  \nabla g(u)=g'(u)\nabla u,
$
we have
\begin{align}
 \|\nabla \mu(\cdot,t_k)\|&\leq\|\nabla \Delta \Phi(\cdot,t_k)\|+\|f'(\Phi(\cdot,t_k))\|_{L^\infty}\|\nabla \Phi(\cdot,t_k)\|\nonumber\\
 &\leq (1+\bar{g}(\mathcal{R}))\mathcal{R}:=\bar{C}_5,\label{muinfty}\\
 \|\nabla f(B\phi^{k-1})\|&\leq\|f'(B\phi^{k-1})\|_{L^\infty}\|\nabla B\phi^{k-1}\|\leq \bar{g}(M)(1+2r_{\max})M.\nonumber
\end{align}

Thus, one has
\begin{equation*}
\|\nabla\mu(\cdot,t_k)-\nabla\Delta\bar\phi^k+\nabla f(B\phi^{k-1})\|
\leq \|\nabla\mu(\cdot,t_k)\|+\|\nabla\Delta\bar e^k\|+\|\nabla\Delta \Phi(\cdot,t_k)\|+\|\nabla f(B\phi^{k-1})\|=\bar{C}_6+\|\nabla\Delta\bar e^k\|,
\end{equation*}
where $\bar{C}_6=(2+\bar{g}(\mathcal{R}))\mathcal{R}+\bar{g}(M)(1+2r_{\max})M$.
Consequently,
\begin{align}\label{e481}
Q_2^k&=\Big( -\nabla\Delta\bar\phi^k+\nabla f(B\phi^{k-1})+\nabla\mu(t_k),-\nabla\Delta\bar\phi^k+\nabla f(B\phi^{k-1})-\nabla\mu(t_k)\Big)\nonumber\\
&\leq \|\nabla\mu(t_k)-\nabla\Delta\bar\phi^k+\nabla f(B\phi^{k-1})\|
\|-\nabla\Delta\bar e^k+\nabla E_f^k-\nabla R_f^k\|\nonumber\\
&\leq (\bar{C}_6+\|\nabla\Delta \bar e^k\|)(\|\nabla\Delta \bar e^k\|+\|\nabla E_f^k-\nabla R_f^k\|).
\end{align}

With the help of the Cauchy-Schwarz inequality, we have
\begin{equation}\label{ca_sc}
\sum\limits_{k=1}^n\sqrt{\tau_k}\|\nabla\Delta \bar e^k\|\cdot \sqrt{\tau_k}
\leq (\sum\limits_{k=1}^n\tau_k\|\nabla\Delta \bar e^k\|^2)^{\frac{1}{2}}(\sum\limits_{k=1}^n\tau_k)^{\frac{1}{2}}.
\end{equation}

Combining the estimates \eqref{e481} and \eqref{ca_sc}, we have
\begin{align}
\sum\limits_{k=1}^n\tau_k Q_2^k\leq&\sum\limits_{k=1}^n\tau_k(\bar{C}_6+\|\nabla\Delta \bar e^k\|)(\|\nabla\Delta \bar e^k\|+\|\nabla E_f^k-\nabla R_f^k\|)\nonumber\\
\leq&\sum\limits_{k=1}^n\tau_k\|\nabla\Delta \bar e^k\|^2+\bar{C}_6t_n\max_{1\leq k\leq n}\|\nabla E_f^k-\nabla R_f^k\|\nonumber\\
&+(\sum\limits_{k=1}^n\tau_k\|\nabla\Delta \bar e^k\|^2)^{\frac{1}{2}}\sqrt{t_n}(\bar{C}_6+\max_{1\leq k\leq n}\|\nabla E_f^k-\nabla R_f^k\|)\leq C_5(\tau+N^{3-m}),
\end{align}
where the last step follows from \eqref{result2}-\eqref{Ej} and  $C_5=C_2+\bar{C}_6t_nC_3+\sqrt{C_2T}(\bar{C}_6+C_3)$.
For the last term of \eqref{err_gamma}, we can derive
\begin{align}\label{Gam_1}
\sum_{k=1}^{n}\tau_k R_\gamma^k&=\sum_{k=1}^{n}[\tau_{k}\partial_t\Gamma(t_k)-\Gamma(t_k)+\Gamma(t_{k-1})]=-\sum_{k=1}^n \int^{t_k}_{t_{k-1}}(t_k-s)\partial_{tt}\Gamma(s){\rm{d}s}\nonumber\\
&\leq \tau\int_0^T|\partial_{tt}\Gamma(s)|{\rm{d}s}\leq 2T(\bar{g}(\mathcal{R})\mathcal{R}+1)\mathcal{R}\tau:=\bar{C}_7\tau,
\end{align}
where the last inequality follows from the fact
\begin{align*}
\partial_{tt}\Gamma(t)&=\partial_t(-\Delta \Phi(\cdot,t)+f(\Phi(\cdot,t))),\partial_t\Phi(\cdot,t))\\
&=(-\Delta \partial_t\Phi(\cdot,t)+f'(\Phi(\cdot,t))\partial_t\Phi(\cdot,t),\partial_t\Phi(\cdot,t))+(-\Delta \Phi(\cdot,t)+f(\Phi(\cdot,t)),\partial_{tt}\Phi(\cdot,t))\\
&\leq(\|\Delta \partial_t\Phi(\cdot,t)\|+\|f'(\Phi(\cdot,t))\partial_t\Phi(\cdot,t)\|)\|\partial_t\Phi(\cdot,t)\|+(\|\Delta \Phi(\cdot,t)\|+\|f(\Phi(\cdot,t))\|)\|\partial_{tt}\Phi(\cdot,t)\|\\
&\leq 2(\bar{g}(\mathcal{R})\mathcal{R}+1)\mathcal{R}.
\end{align*}
Combining the estimates \eqref{Q_1}-\eqref{Gam_1} and \eqref{err_gamma}, we have
$$s^n\leq \bar{C}_5^2\sum_{k=1}^n\tau_{k}s^{k} +[C_4t_n+\gamma^0C_5+\bar{C}_7](\tau+N^{3-m})+s^0.$$
With $\tau < 1/(\bar{C}_5^2+1)$ in \eqref{Limit_tau} and initial condition $\gamma^0=\Gamma(t_0)$, it is a consequence of Lemma \ref{Gronwall} that
\begin{equation}\label{result3}
|s^n|\leq C_6 (\tau+N^{3-m}).
\end{equation}
Here the constant $C_6$ can be defined as
$$C_6=\exp((1+\bar{C}_5^2)\bar{C}_5^2T)(1+\bar{C}_5^2)(C_4T+\gamma^0C_5+\bar{C}_7).$$

\textbf{Step 5: Estimate of} $|1-\xi^n|$. Combining \eqref{scheme3}, \eqref{E_E} and \eqref{result3} yields
\begin{align}
|1-\xi^n|&=|\frac{E(\bar\phi^n)+1-\gamma^n}{E(\bar\phi^n)+1}|
=|\frac{E(\bar\phi^n)-E(\Phi(\cdot,t_n))+\Gamma^n-\gamma^n}{E(\bar\phi^n)+1}|\nonumber\\
&\leq C_4(\tau+N^{3-m})+s^n\leq (C_4+C_6)(\tau+N^{3-m}).
\end{align}
We set $c_0=C_4+C_6$, then the mathematical induction for \eqref{conclution} is completed.
In fact, $c_0$ is a constant independent of $n, \tau, N$, which plays an important role in mathematical induction.

Now we substitute $c_0=C_4+C_6$ into the second row of \eqref{result12} and deduce
\begin{equation}
\|\nabla \bar e^n\|\leq C_1(\tau^2+N^{1-m}+(C_4+C_6)^2(\tau+N^{3-m})^2)+C_p\mathcal{R}N^{1-m}\leq C(\tau^2+N^{3-m}).\label{barnabla_e}
\end{equation}
Similarly, from the third inequality of \eqref{result22}, we have
\begin{equation}
\sum_{k=1}^n\tau_k\|\nabla \Delta \bar e^k\|^2\leq C(\tau^4+N^{2(3-m)}).\label{barnablaDelta_e}
\end{equation}
Noticing $e^n=\bar e^n+(1-\eta^n)\bar \phi^n$ and \eqref{tri2},  $\|\nabla e^n\|$ and $\sum\limits_{k=1}^n\tau_k\|\nabla\Delta e^k\|^2$ can be estimate as \eqref{re1}-\eqref{re2}  by the triangle inequality and \eqref{barnabla_e},\eqref{barnablaDelta_e}.
Besides, \eqref{re3} follows from \eqref{result3}.
The proof is completed.
\end{proof}

\begin{remark}\label{remark2}
We point out that the inequality zoom in \eqref{est_E} plays an important role in achieving the optimal  second-order error estimate.
In fact, if the Gr\"onwall's inequality is used directly, the Young's inequality leads to $\sum\limits_{k=1}^n\tau_k\|\nabla \Delta R_{f}^j+\nabla R_{t}^j\|^2=O(\tau^3)$ inevitably as discussed in \cite{mayuheng2}. This ultimately leads to the order reduction as follows
\begin{equation*}
\|\nabla \bar e_N^n\|\leq C_0\left(\tau^{3/2}+N^{1-m}+c_0^2(\tau+N^{3-m})^2\right),
\end{equation*}
instead of \eqref{e4660}.
On the contrary, by using the inequality zoom, the term $\|\nabla\bar e_N^{n_0}\|$ can be eliminated in \eqref{err}, which avoids the use of Young's inequality. In this case, the Gr\"onwall's inequality obtains the  second-order sharp estimate of $\|\nabla \bar e_N^n\|$.
\end{remark}

%
%
\section{Numerical examples}
We now present 2D and 3D numerical examples to demonstrate the accuracy, energy stability, and efficiency of the proposed scheme \eqref{scheme}. In simulations, we set the computational domain $\Omega=[0,2\pi]^d (d=2,3)$ with periodic boundary conditions.  We adopt the Fourier-spectral method to discretize the space with $N$ Fourier modes for each directions.
Unless otherwise specified, we take $N=128$.

\subsection{Accuracy test of time}
We first investigate the temporal convergence rate for scheme \eqref{scheme} in 2D.
Here we take $\varepsilon=0.2$ and the initial value as
\begin{equation*}
  \Phi^0(\bm{x})=-\tanh\big(\frac{\sqrt{(x-\pi)^2+(y-\pi)^2}-1.5}{4\varepsilon}\big).
\end{equation*}
Since the exact solutions of the system are not known, we compute the reference solution by scheme \eqref{scheme} with a tiny time size $\tau=1\times 10^{-7}$.
We generate a random time mesh $\tau_k=T\theta_k/(\sum_{k=1}^{K}\theta_k)$ under condition \Ass{1}, where $T$ is the final time and $\theta_k$  is a random perturbation uniformly distributed in
$(1/4.86,1)$.
Denote $e$ by the error. The corresponding convergence order of time at $T$ is calculated by
 \begin{equation*}
 \text{Order}:= \frac{\log e(K)-\log e(2K) }{ \log \tau(K)-\log \tau(2K)},
 \end{equation*}
where $\tau(K)$ represents the max time step size of $\{\tau_k\}_{k=1}^K$.

By varying the number of random steps $K$, we compute the $H^1$ errors of the phase variable and the absolute error of modified energy between the reference and approximate solutions at time $t = 0.1$.
As shown in Table \ref{table1}, $\|\phi^K-\Phi(\cdot,t_K)\|_{H^1}$  and $|\gamma^K-\Gamma(t_K)|$ have the second-order and first-order convergence accuracy in time, respectively, which coincides with the theoretical analysis in Theorem \ref{Th5}.
\begin{table}[!ht]
\begin{center}
\caption{ The error and convergence order of BDF2 scheme (\ref{scheme}) at $t=0.1$}
\def\temptablewidth{1.0\textwidth}
{\rule{\temptablewidth}{1pt}}
\begin{tabular*}{\temptablewidth}{@{\extracolsep{\fill}}ccccccc}
 $K$ &$\tau$ &$\|\phi^K-\Phi(\cdot,t_K)\|_{H^1}$ & Order&$|\gamma^K-\Gamma(t_K)|$ & Order& max $r_k$\\
\hline
400&4.0731e-04&5.6229e-04&--&4.5841e-01&--&4.501\\
800&2.0643e-04&1.5974e-04&1.85&2.4304e-01&0.93&4.799\\
1600&1.0229e-04&3.6817e-05&2.09&1.1966e-01&1.01&4.731\\
3200&5.1749e-05&9.7311e-06&1.95&6.0611e-02&1.00&4.774\\
\end{tabular*}
{\rule{\temptablewidth}{1pt}}
\label{table1}
\end{center}
\end{table}
\subsection{Coalescence of two kissing bubbles in 2D}
We now consider the coalescence of two kissing bubbles by taking $\varepsilon^2=0.1$ and the initial value as
\begin{equation*}
  \Phi^0(\bm{x})=\sum_{i=1}^2\tanh\big(r_i-\frac{\sqrt{(x-x_i)^2+(y-y_i)^2}}{4\varepsilon^2}\big),
\end{equation*}
where
$x_1=\pi-1,y_1=\pi-\pi,r_1=1,x_2=\pi+1,y_2=\pi,r_1=1.$

 The adaptive time-stepping strategy \cite{fan,qiao_adaptive_2011} was adopted by
  \begin{equation}
  \tau_{n+1} =\max \left\{r_{\max}\tau_{n},\tau_{\min},\frac{\tau_{\max}}{\sqrt{1+\alpha|\delta_t E(t_n)|^2}}\right\},
  \label{Adaptive}
  \end{equation}
where $\delta_t E(t_n)$ presents the discrete temporal derivative of the energy defined by
 $\delta_t E(t_n) = (\gamma^n-\gamma^{n-1})/\tau_n$.
 Here $\tau_{\min}$ and $\tau_{\max}$ are the minimum and  maximum time steps respectively and $\alpha$ is a tunable parameter related to the  level of the adaptivity.
Here we take $\tau_{\min}=10^{-4}$, $\tau_{\max}=7\times10^{-3}$ and $\alpha=0.01$.

\begin{figure}[!ht]
\centering
\subfigure[]{
\includegraphics[width=5cm]{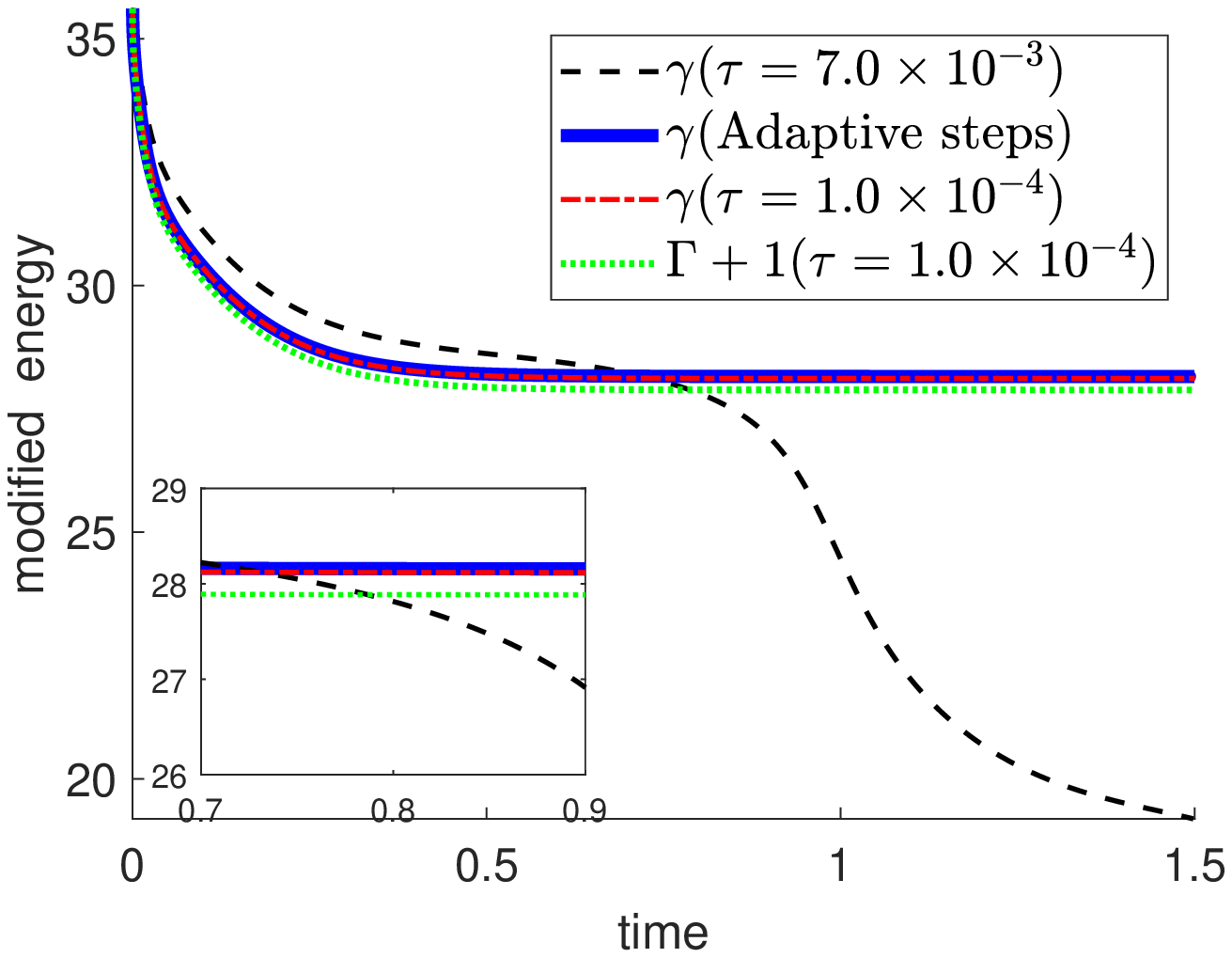}
\label{fig_energy}}
\subfigure[]{
\includegraphics[width=5cm]{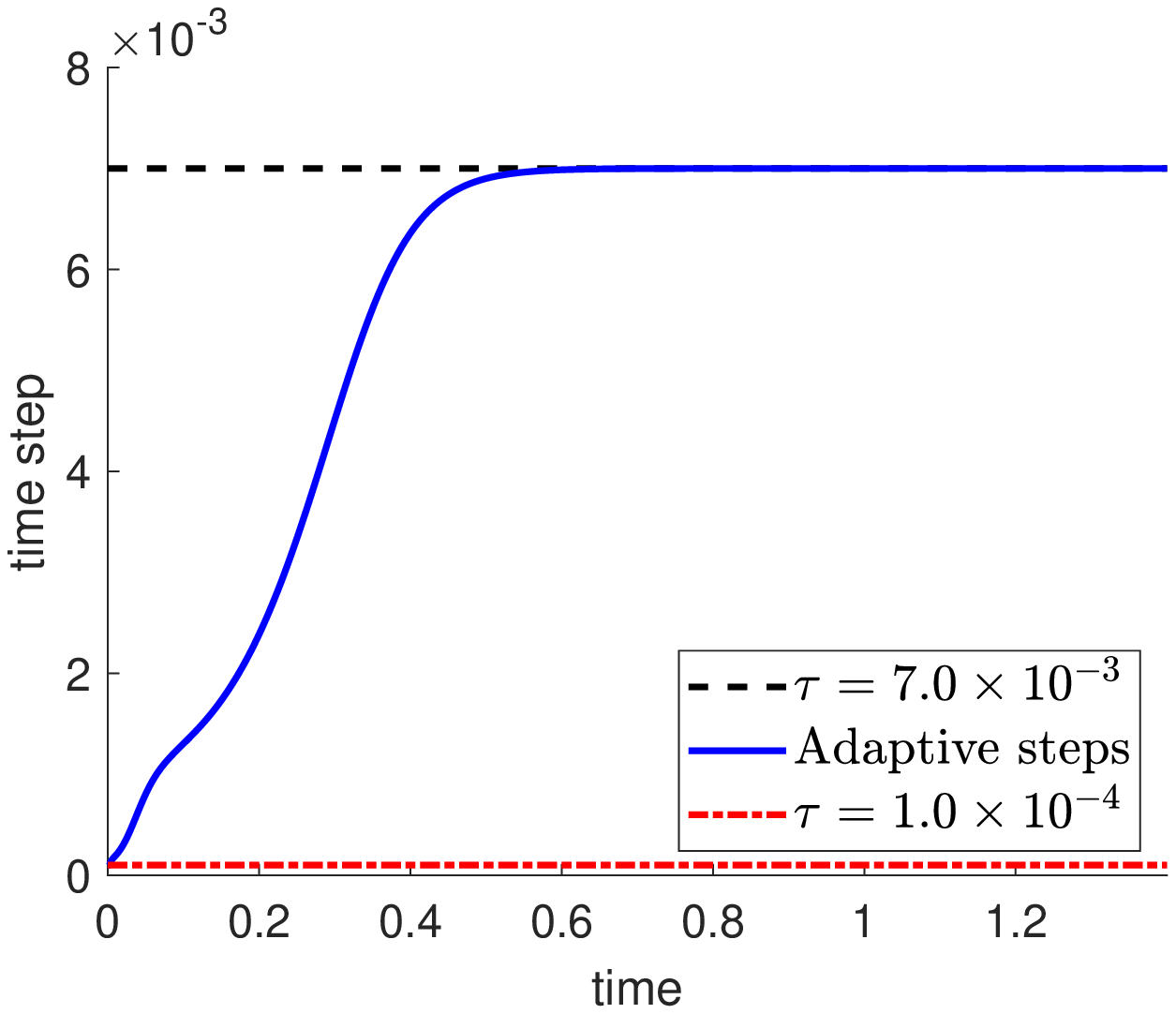}
\label{timestep}}
\subfigure[]{
\includegraphics[width=5cm]{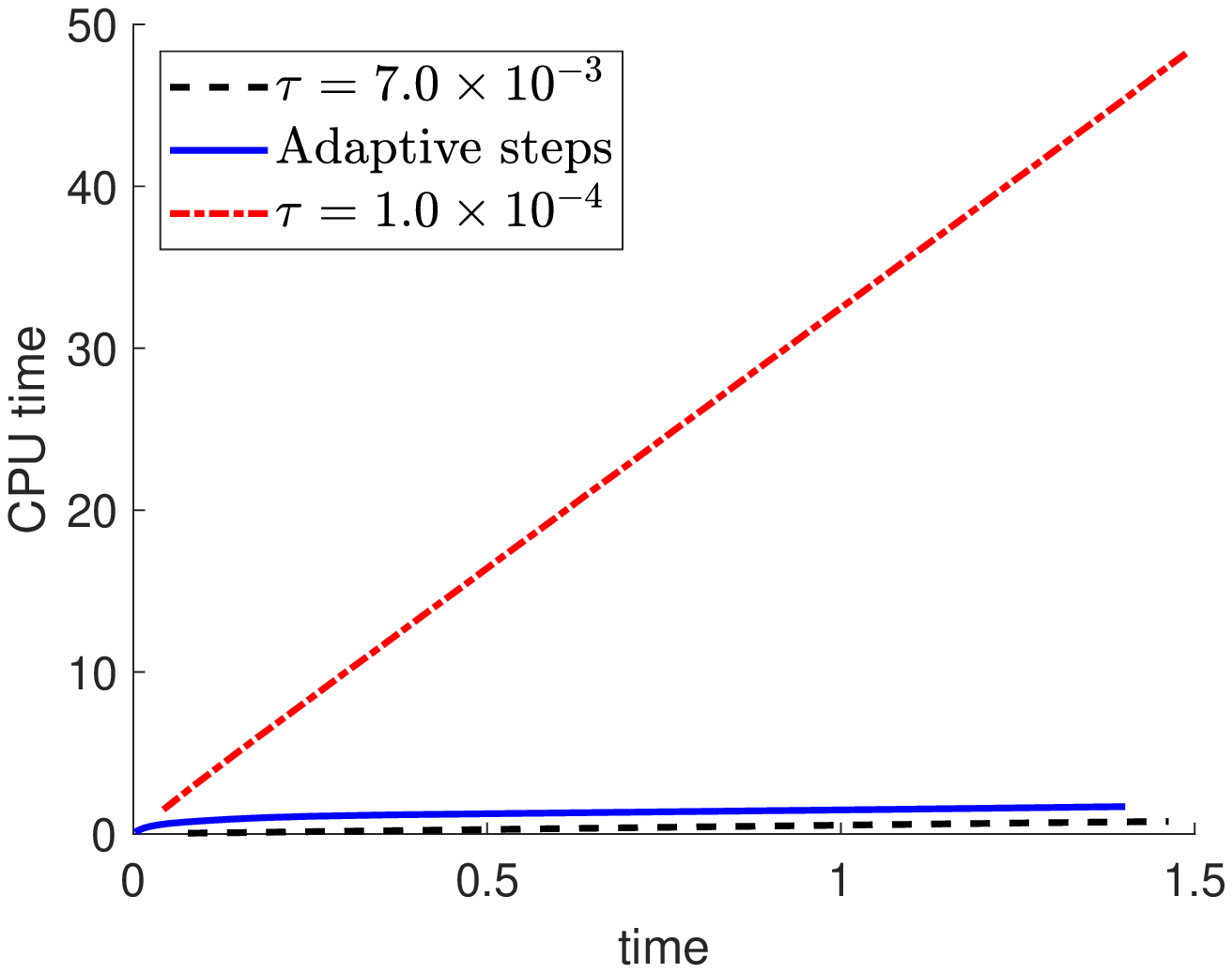}
\label{CPU_2D}}
\caption{(a): Evolution of modified discrete energy; (b): adaptive time-step size; (c): CPU time contrast}
\end{figure}

We compare the evolution of the modified energy $\gamma^n$ under different time steps, and select the original energy with $\tau=10^{-4}$ as the reference solution.
Figure \ref{fig_energy} shows that $\gamma^n$ converges to the wrong solution when the time step is chosen too large, such as $\tau=7\times10^{-3}$.
There is no difference in the evolution of $\gamma^n$ between the adaptive step and $\tau=10^{-4}$.
Meanwhile, the modified energy $\gamma^n$ is decreasing all the time in keeping with the evolution of $\Gamma(t_n)$.
This verifies that the given adaptive scheme preserves the energy dissipation law without sacrificing accuracy.

We now consider the CPU time for different adaptive and fixed step sizes.
Figure \ref{timestep} plots the adaptive time steps, which shows the adaptive step takes relatively large time steps when energy changes slowly.
As shown in Figure \ref{CPU_2D}, adaptive strategy significantly reduces the CPU time compared with fixed step $\tau=10^{-4}$. In addition, one can see that the CPU times are almost commensurate between the adaptive strategy and the fixed step  $\tau =7\times10^{-3}$, but the solution will be incorrect for $\tau =7\times10^{-3}$ for long time simulations, and the adaptive strategy \eqref{Adaptive} still work well.

\begin{figure}[!ht]
\centering
\subfigure{
\begin{minipage}[b]{0.31\textwidth}
\includegraphics[width=0.9\textwidth]{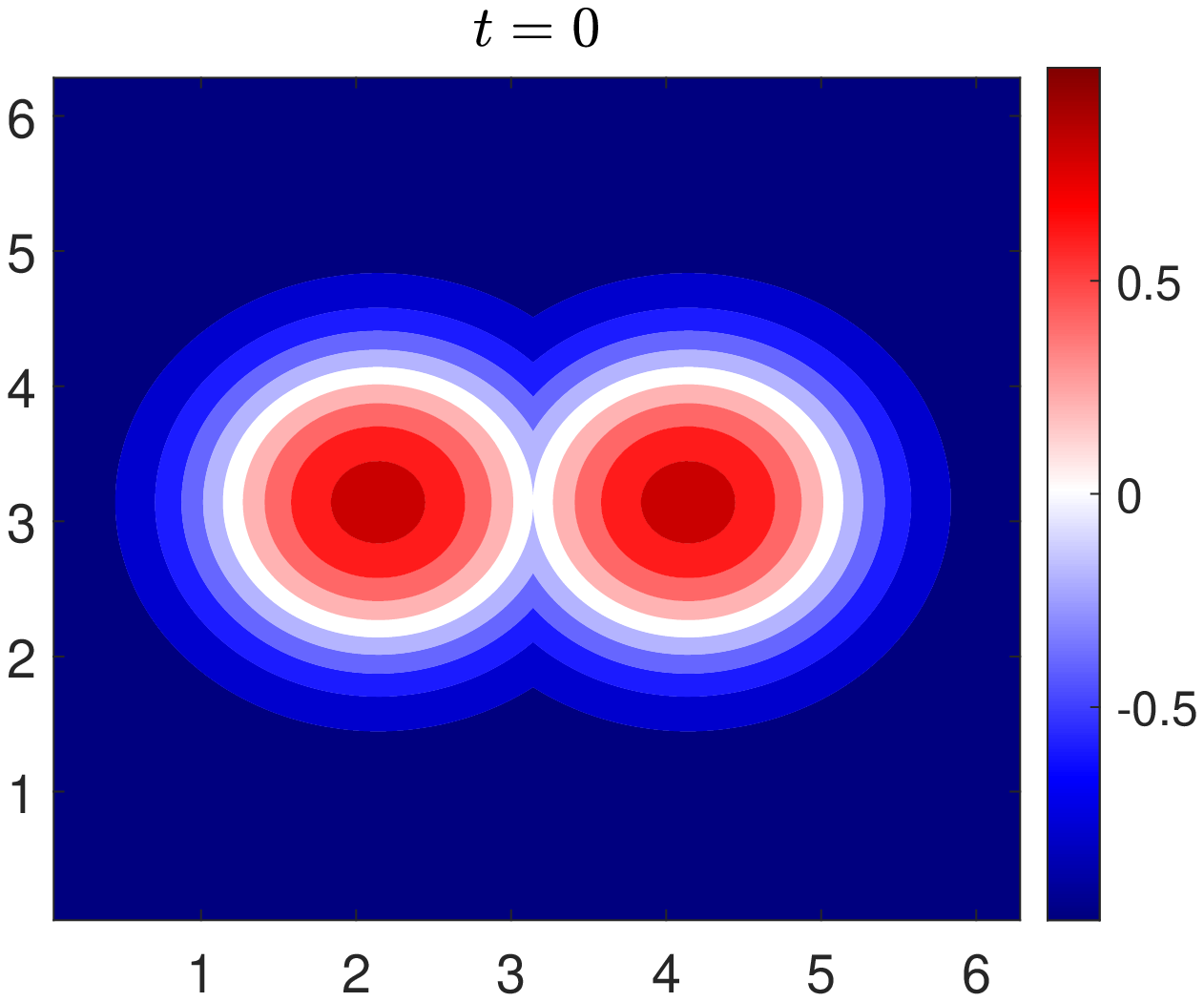}\\
\includegraphics[width=0.9\textwidth]{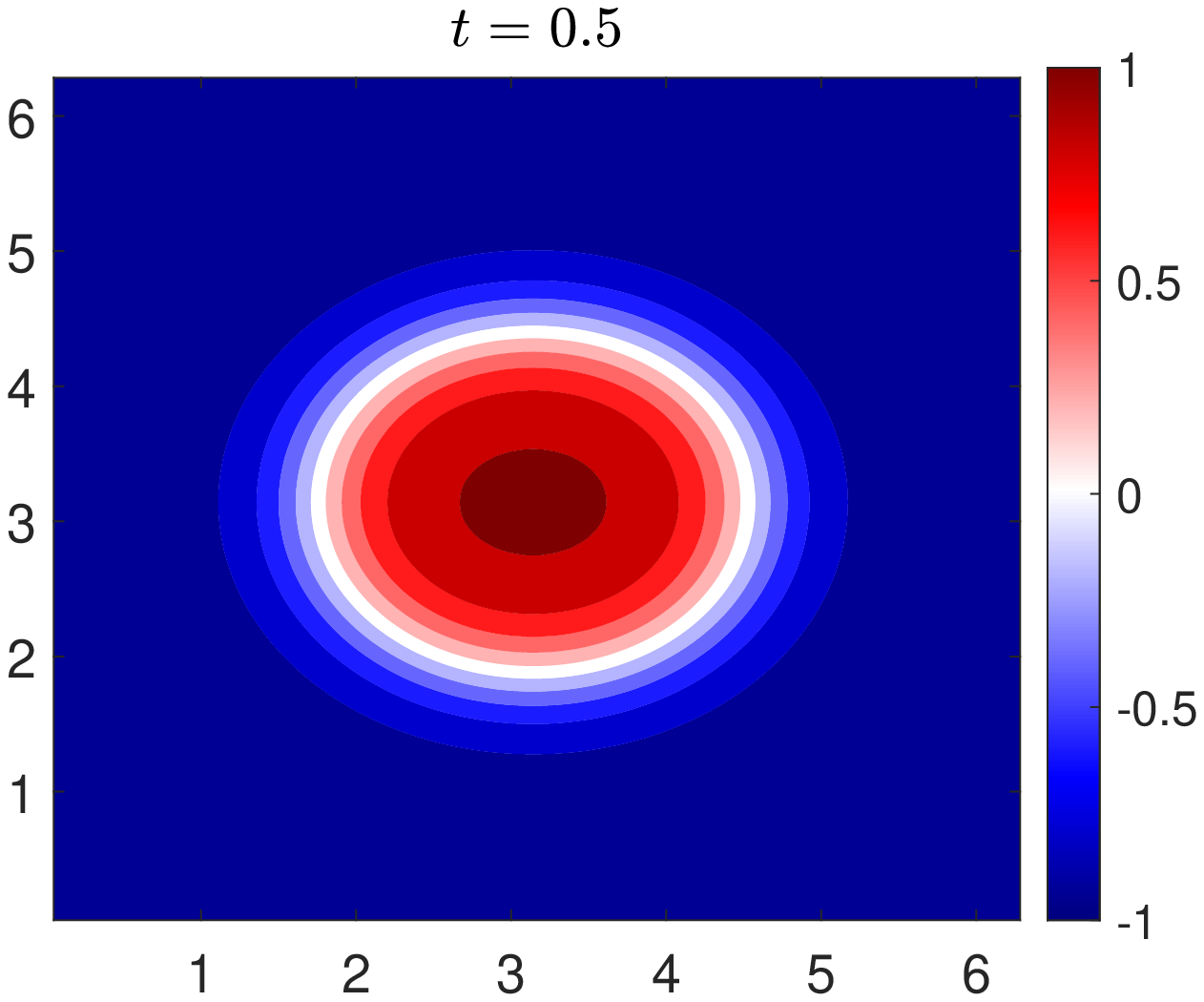}
\end{minipage}
}\hspace{-5mm}
\subfigure{
\begin{minipage}[b]{0.31\textwidth}
\includegraphics[width=0.9\textwidth]{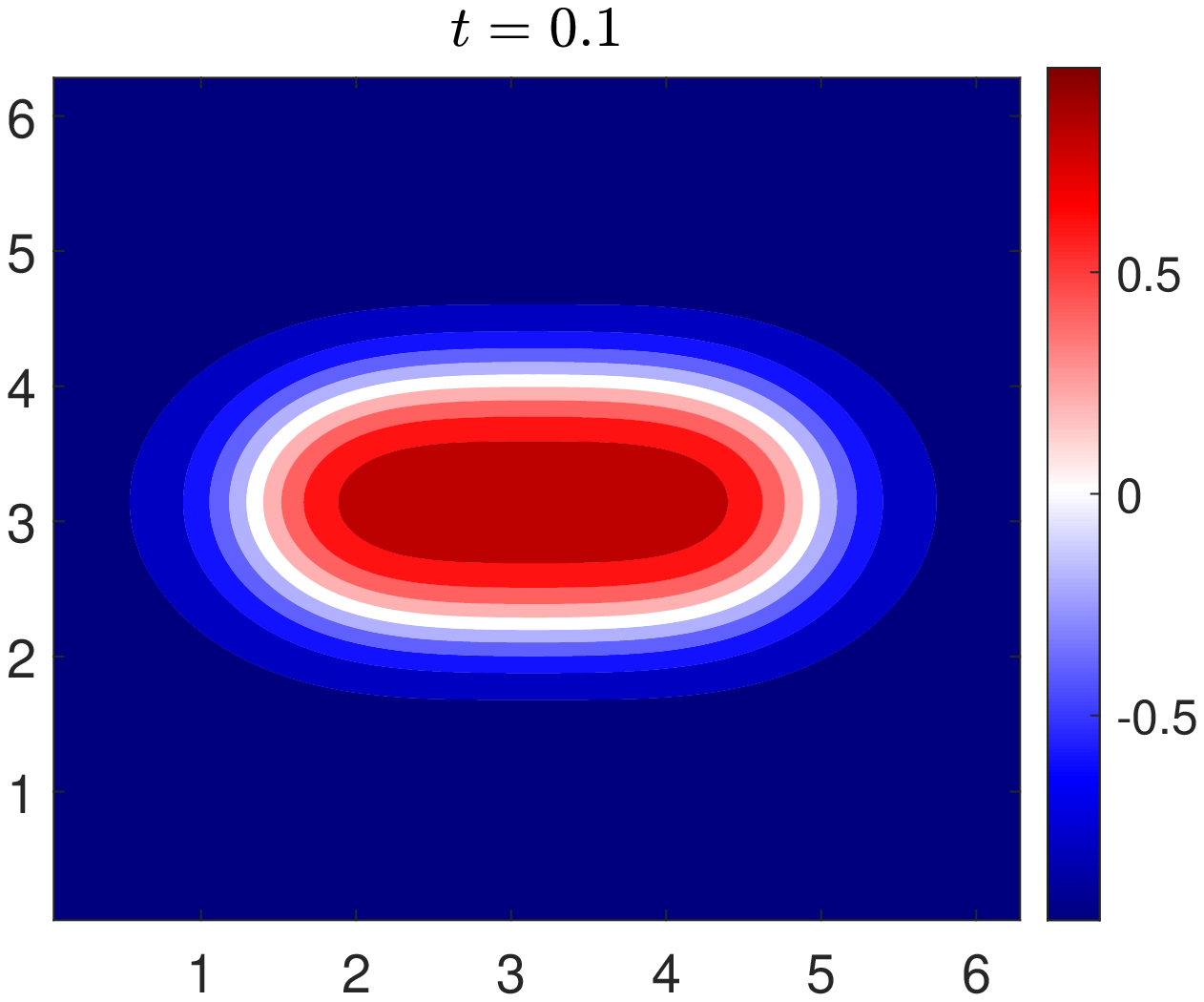} \\
\includegraphics[width=0.9\textwidth]{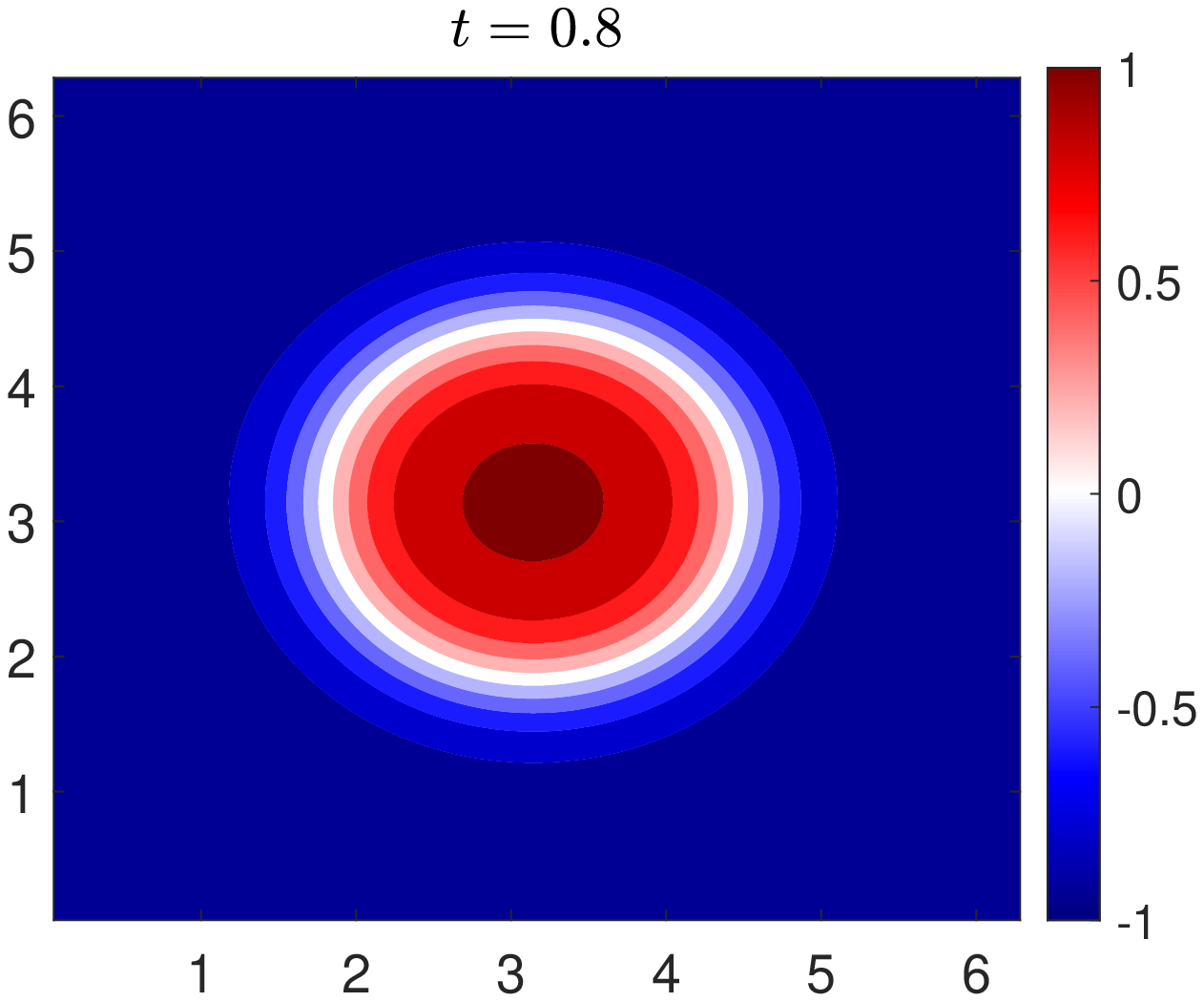}
\end{minipage}
}\hspace{-5mm}
\subfigure{
\begin{minipage}[b]{0.31\textwidth}
\includegraphics[width=0.9\textwidth]{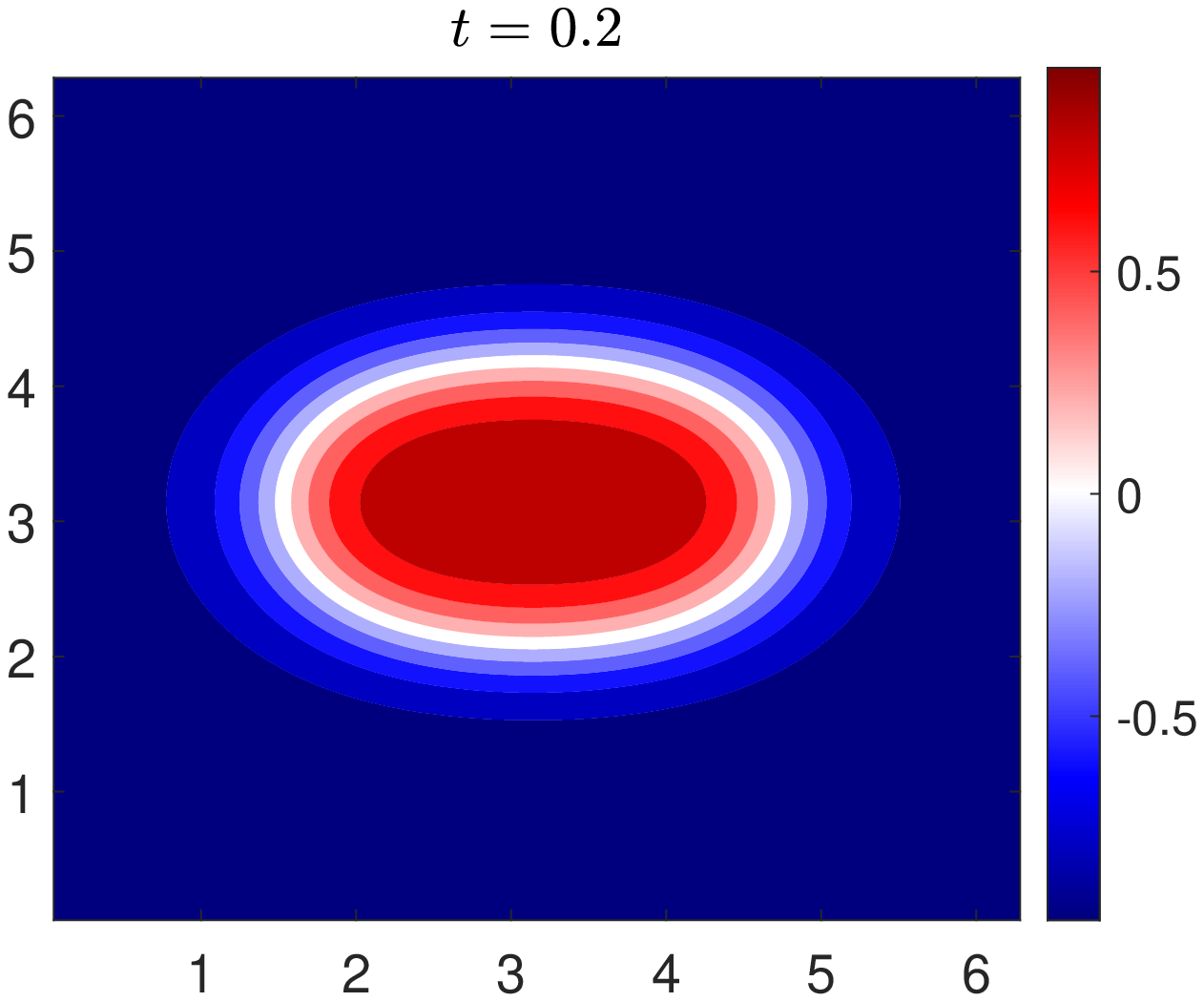} \\
\includegraphics[width=0.9\textwidth]{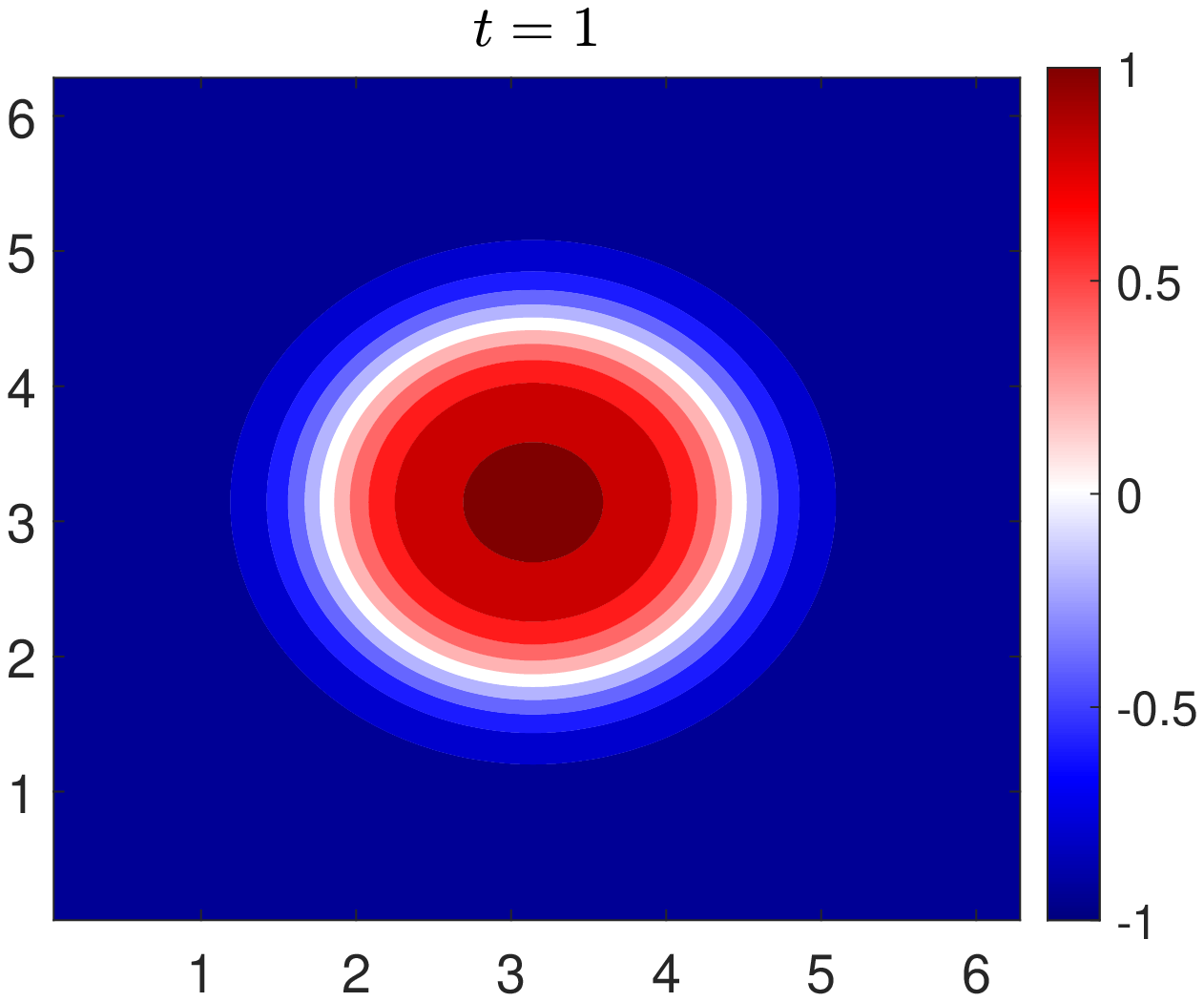}
\end{minipage}
}
\caption{Solution snapshots at t = 0,  0.1, 0.2, 0.5, 0.8, 1.}
\label{snap}
\end{figure}

The phase transition behavior of the density field $\Phi$ is shown in Figure \ref{snap} ,
One can see the two spheres slowly merge and eventually and stabilize into a circular region.
The observed phenomena are consistent with the published results in \cite{zhang_non-iterative_2020}.

\subsection{Evolutions of coarsening process}
We now consider the coarsening dynamics of C-H equation with random initial values in 2D and 3D, respectively.
The random initial value is given as
$$\Phi^0(\bm{x}) = 0.35 + 0.3 {\bf Rand} (\bm{x}) \quad \bm{x}\in (0,2\pi)^d,  d=2,3.$$

\subsubsection{The coarsening process in 2D}
We first investigate the coarsening process in dimension two by taking the parameters as
$\varepsilon=0.3, \tau_{\max}=10^{-4}, \tau_{\min}=10^{-5}, \alpha=0.01.$ The time adaptive strategy is also used by \eqref{Adaptive}.

\begin{figure}[!ht]
\centering
\subfigure[]{
\includegraphics[width=6cm]{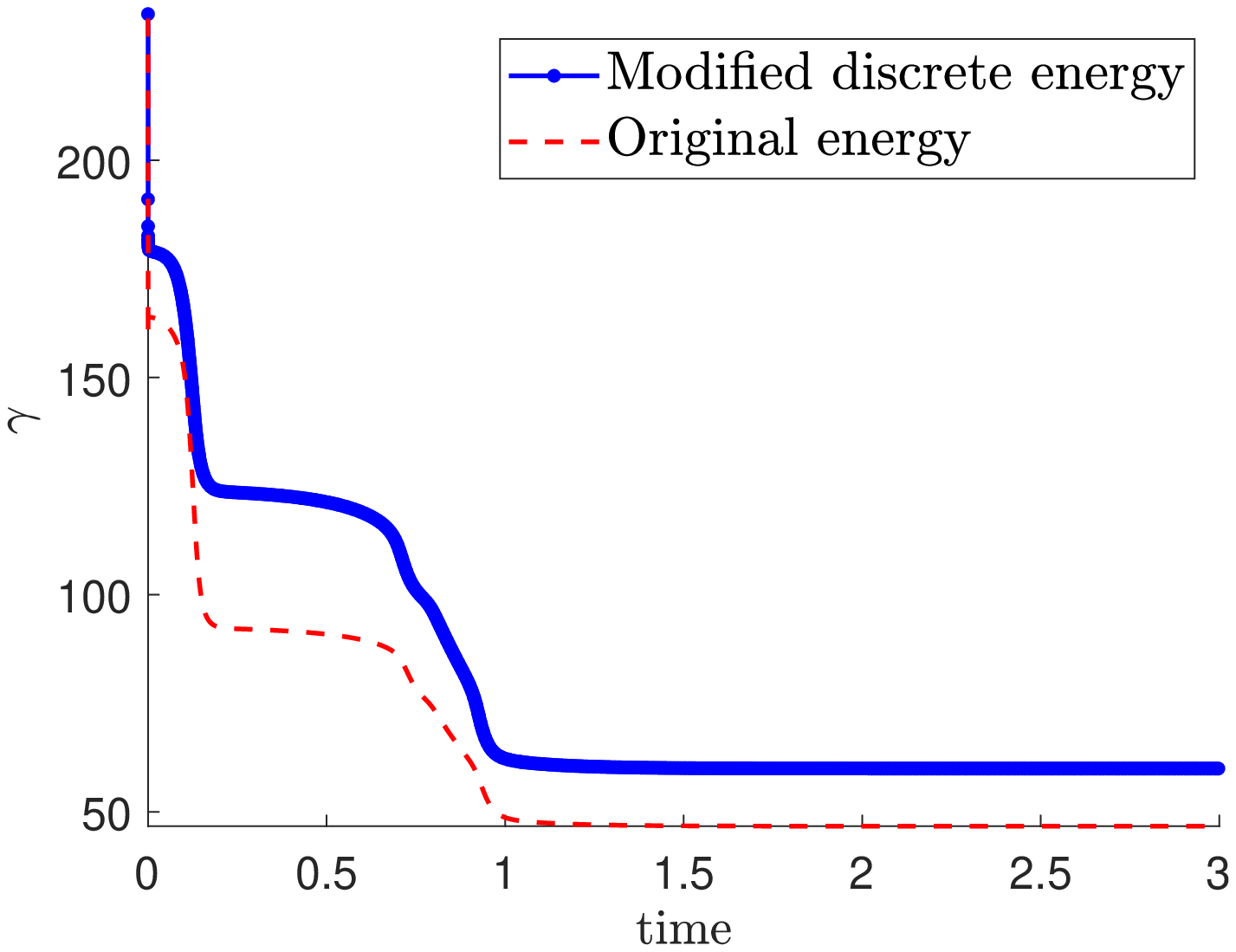}
\label{fig_energy_e2}}
\subfigure[]{
\includegraphics[width=6cm]{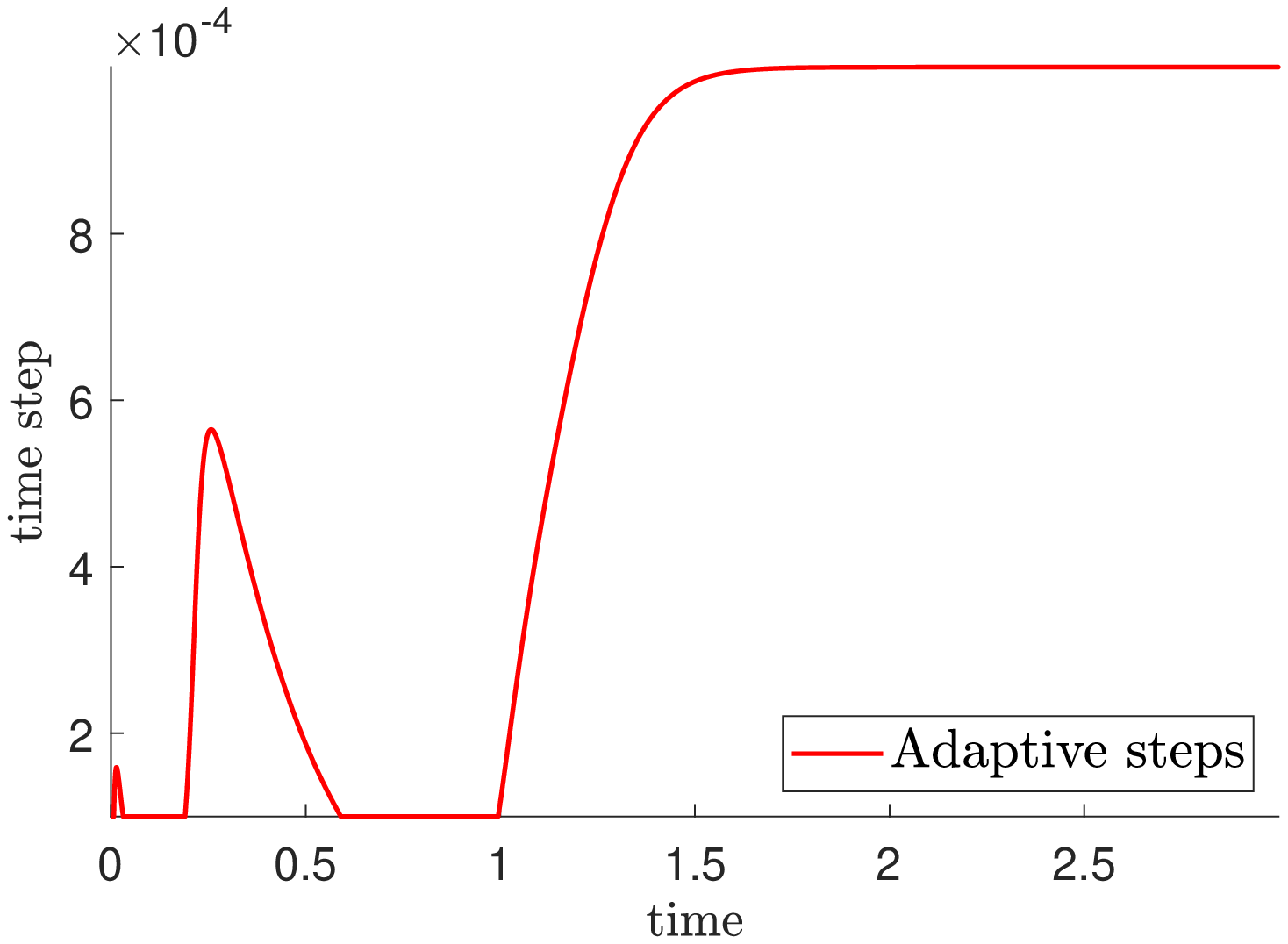}
\label{timestep_e2}}
\caption{Evolutions of modified discrete energy (a) and adaptive time steps (b)}
\end{figure}

\begin{figure}[!ht]
\centering
\subfigure{
\begin{minipage}[b]{0.3\textwidth}
\includegraphics[width=\textwidth]{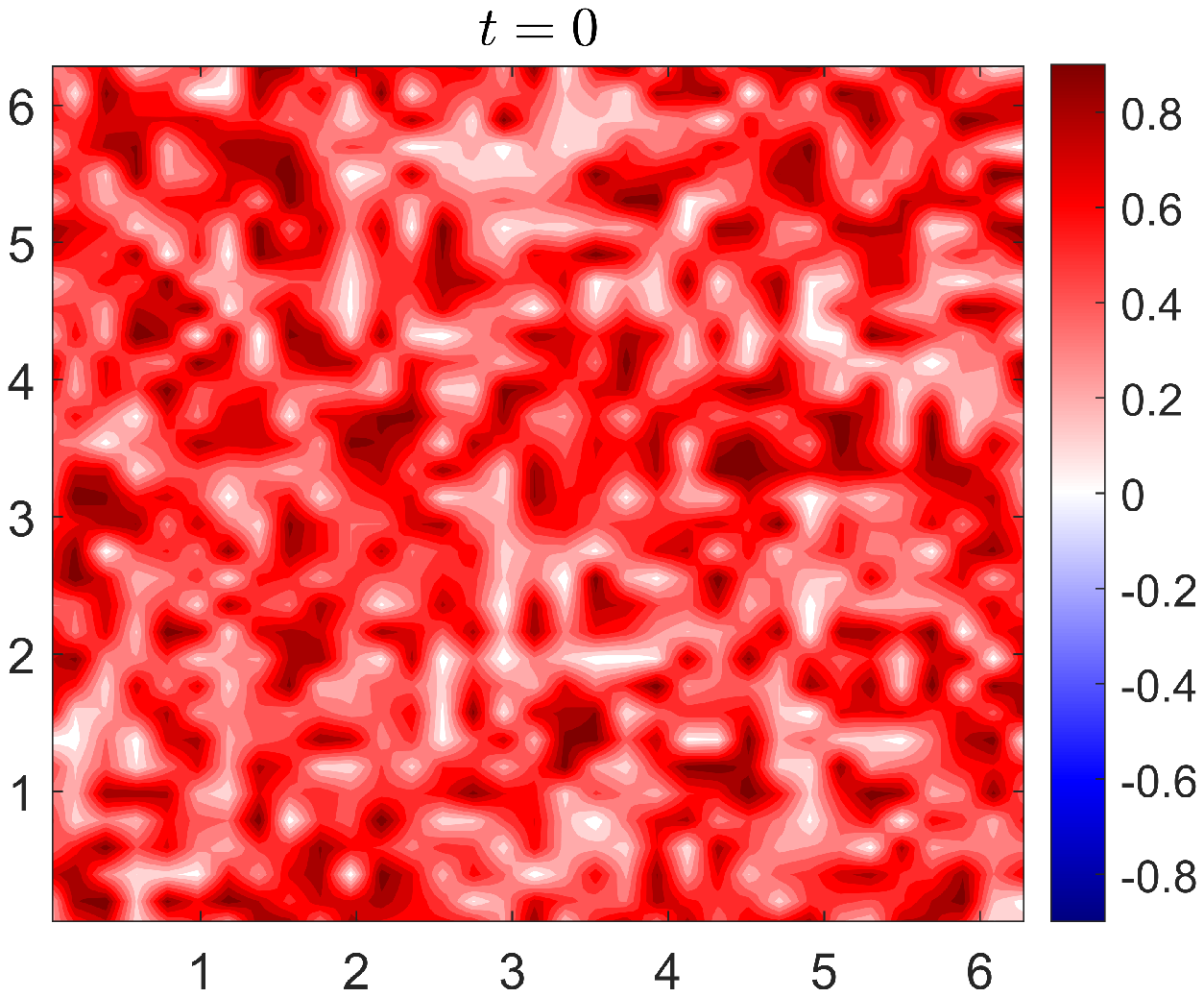}\\
\includegraphics[width=\textwidth]{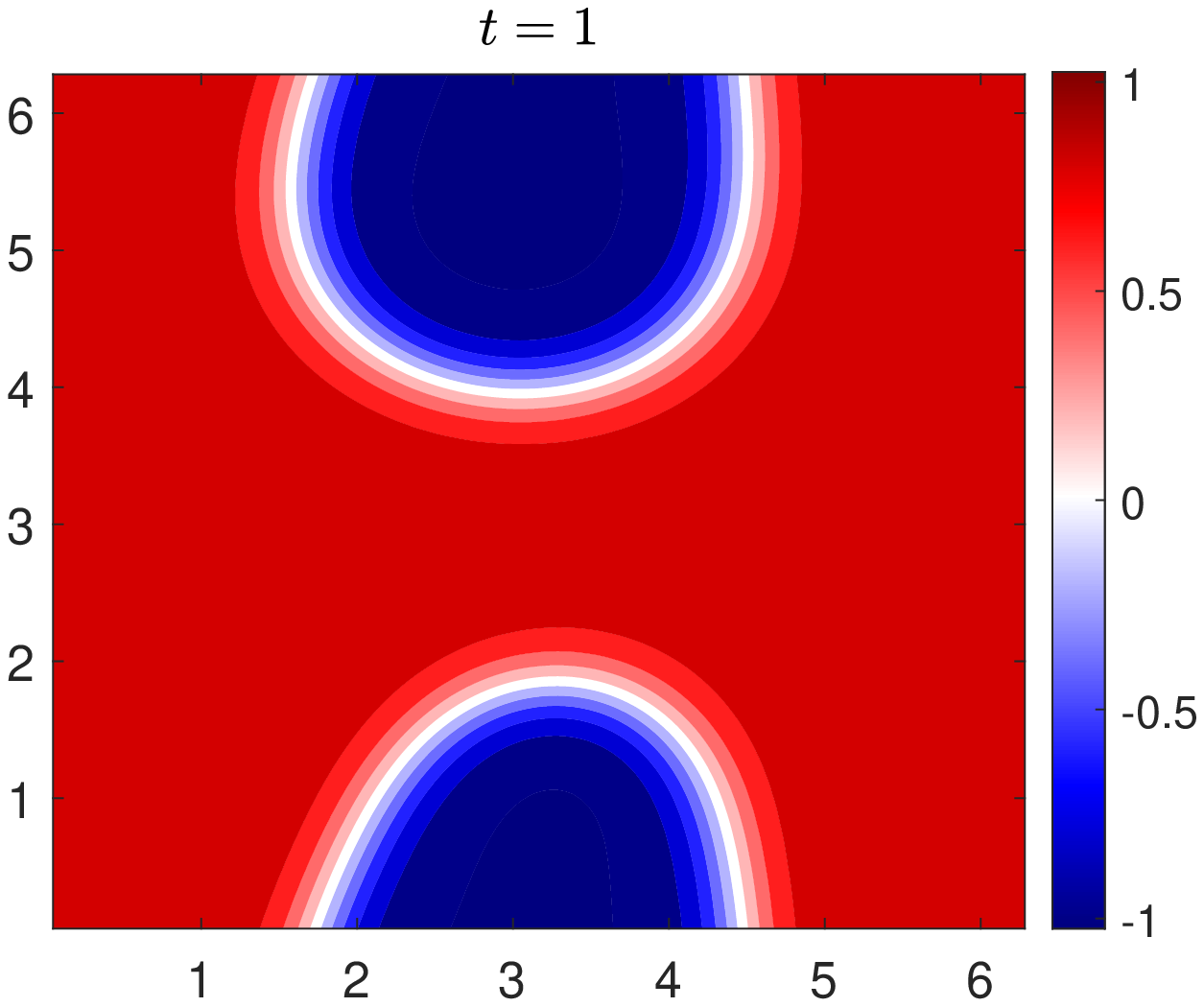}
\end{minipage}
}
\subfigure{
\begin{minipage}[b]{0.3\textwidth}
\includegraphics[width=\textwidth]{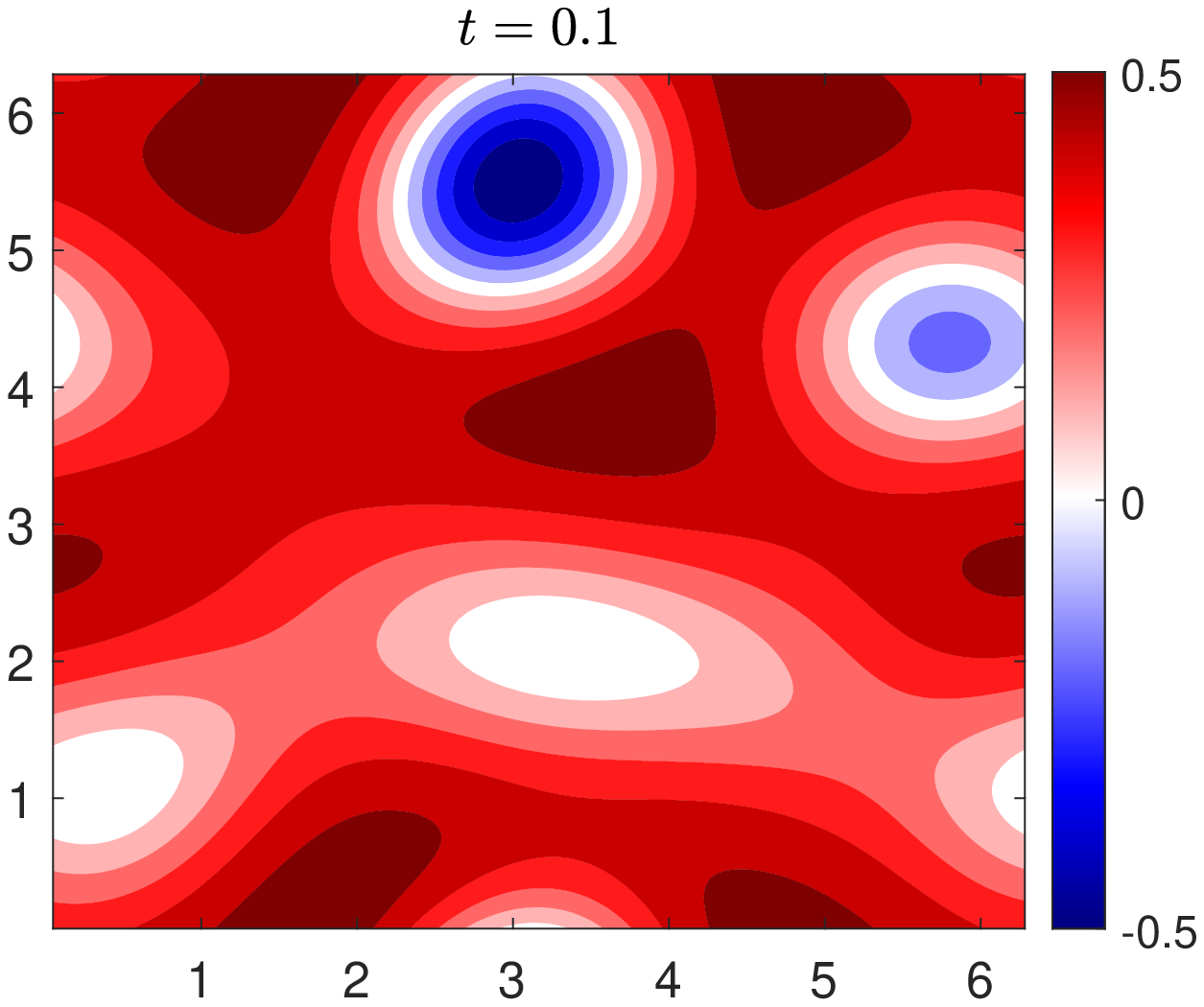} \\
\includegraphics[width=\textwidth]{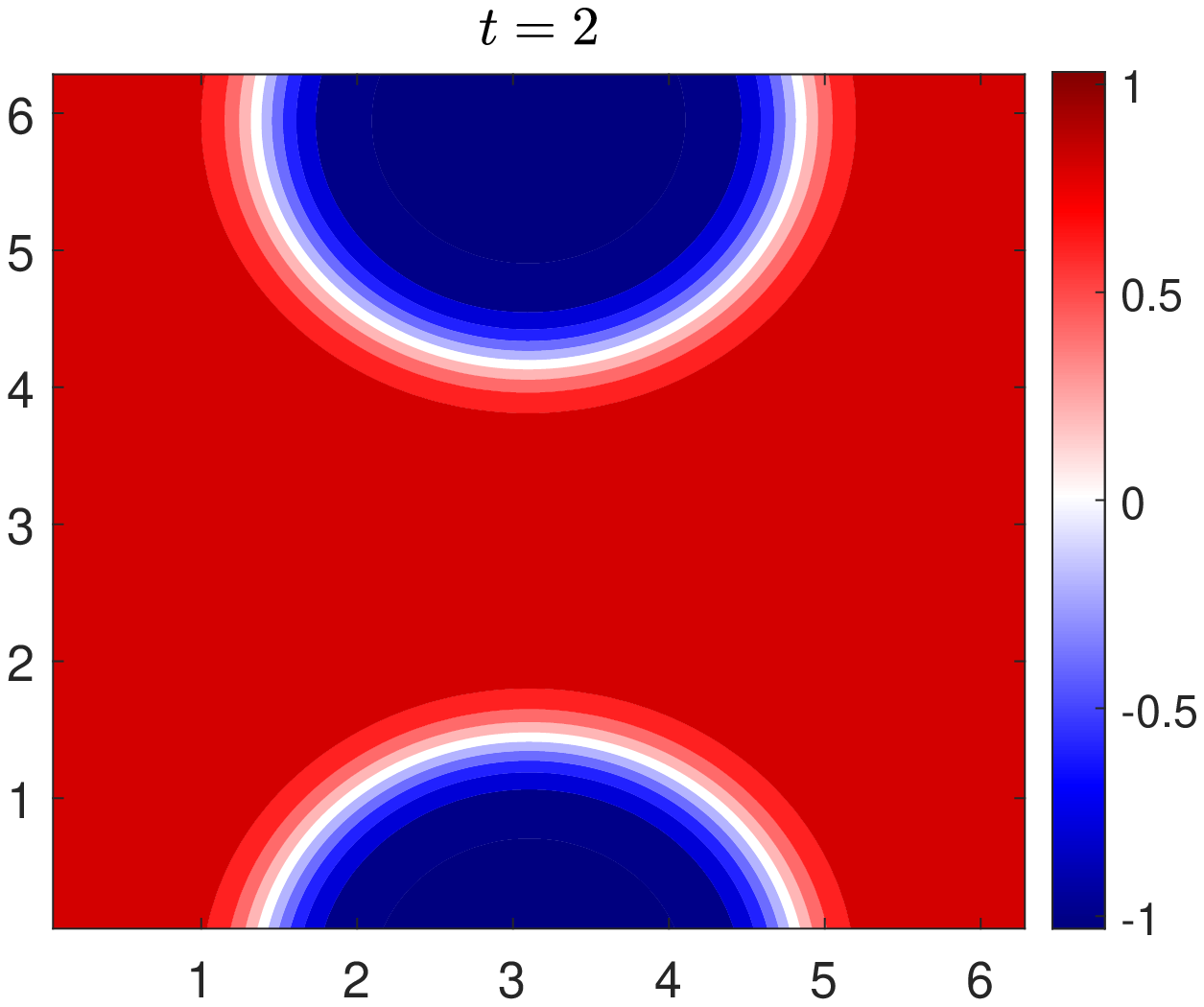}
\end{minipage}
}
\subfigure{
\begin{minipage}[b]{0.3\textwidth}
\includegraphics[width=\textwidth]{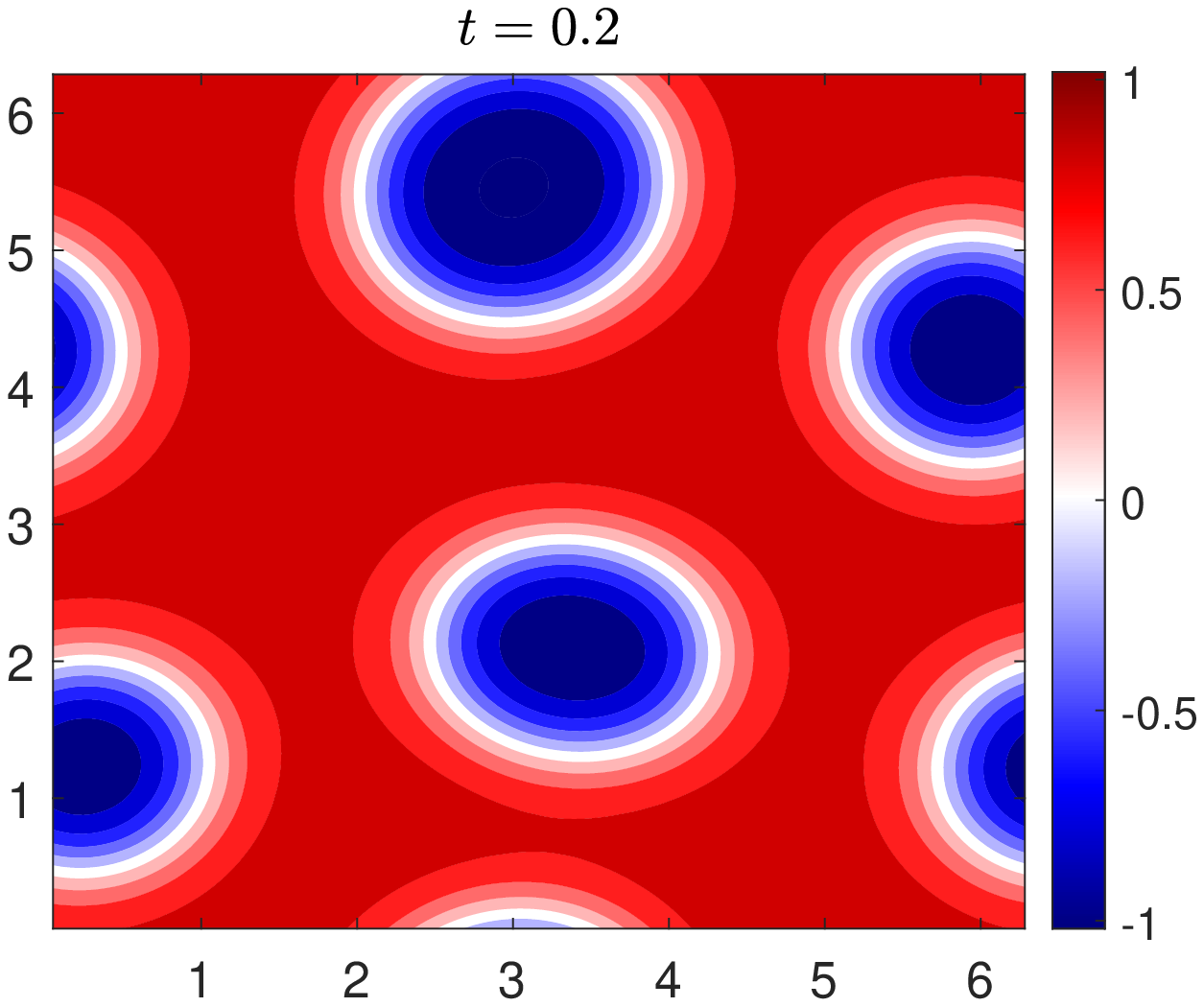} \\
\includegraphics[width=\textwidth]{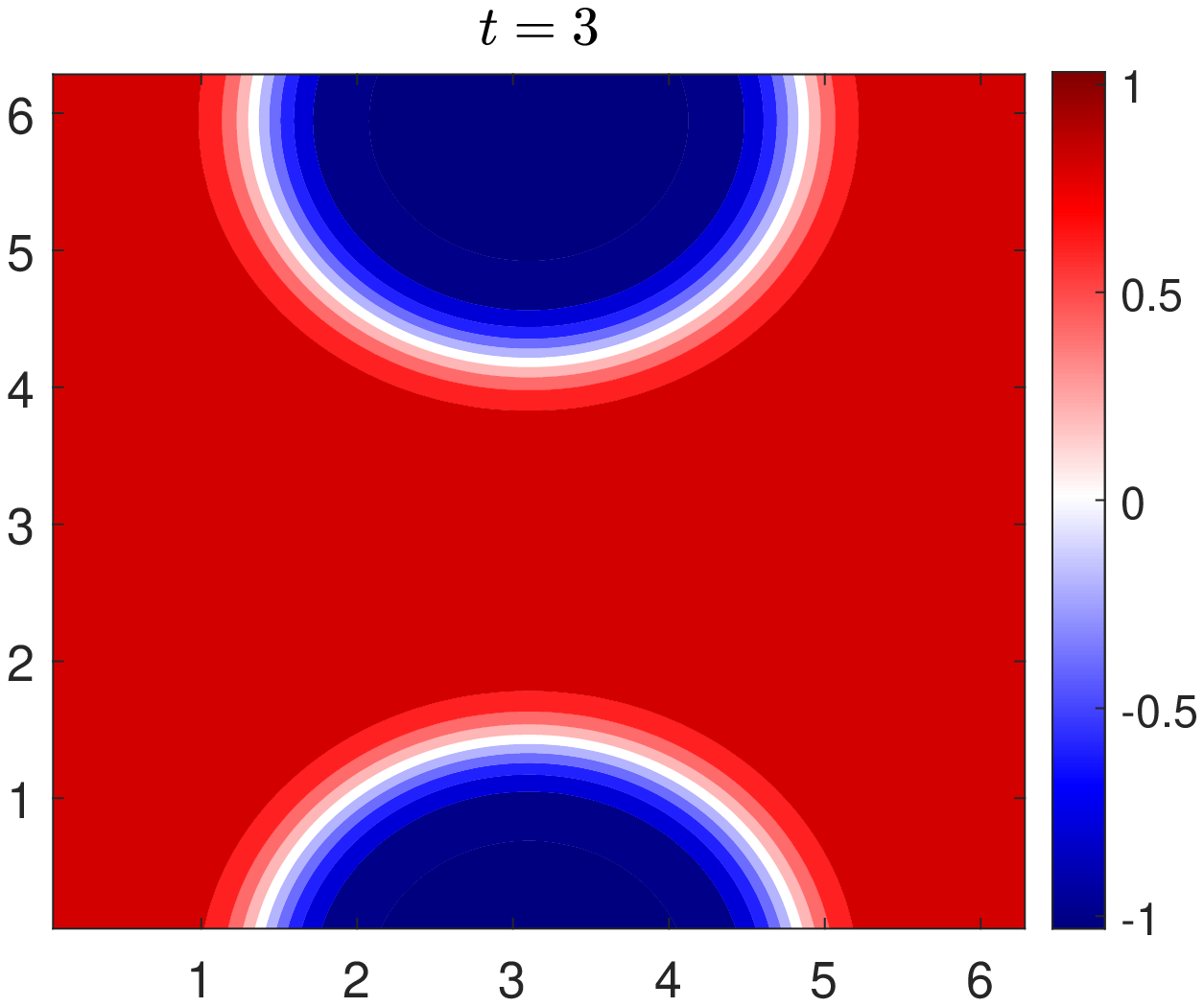}
\end{minipage}
}
\caption{Solution snapshots at t = 0,  0.1, 0.2, 1, 2, 3.}
\label{snap_e2}
\end{figure}

As shown in Figure \ref{fig_energy_e2}, the modified energy is strictly dissipated and closely related to the original energy.
Figure \ref{timestep_e2} shows that  the adaptive step size is relatively small when the energy changes dramatically.
Comparatively, the step size is arger when energy evolves slowly.
Finally, snapshots of the evolution about phase transitions are shown in Figure \ref{snap_e2} which is consistent with our theory.
\subsubsection{The coarsening process in 3D}
We now consider the coarsening process in dimension three.
The parameters are taken as $N=48, \varepsilon=h=2\pi/N, \tau_{\max}=10^{-4}, \tau_{\min}=4\times 10^{-5}, \alpha=1$.

\begin{figure}[!ht]
\centering
\subfigure{
\begin{minipage}[b]{0.31\textwidth}
\includegraphics[width=1.3\textwidth]{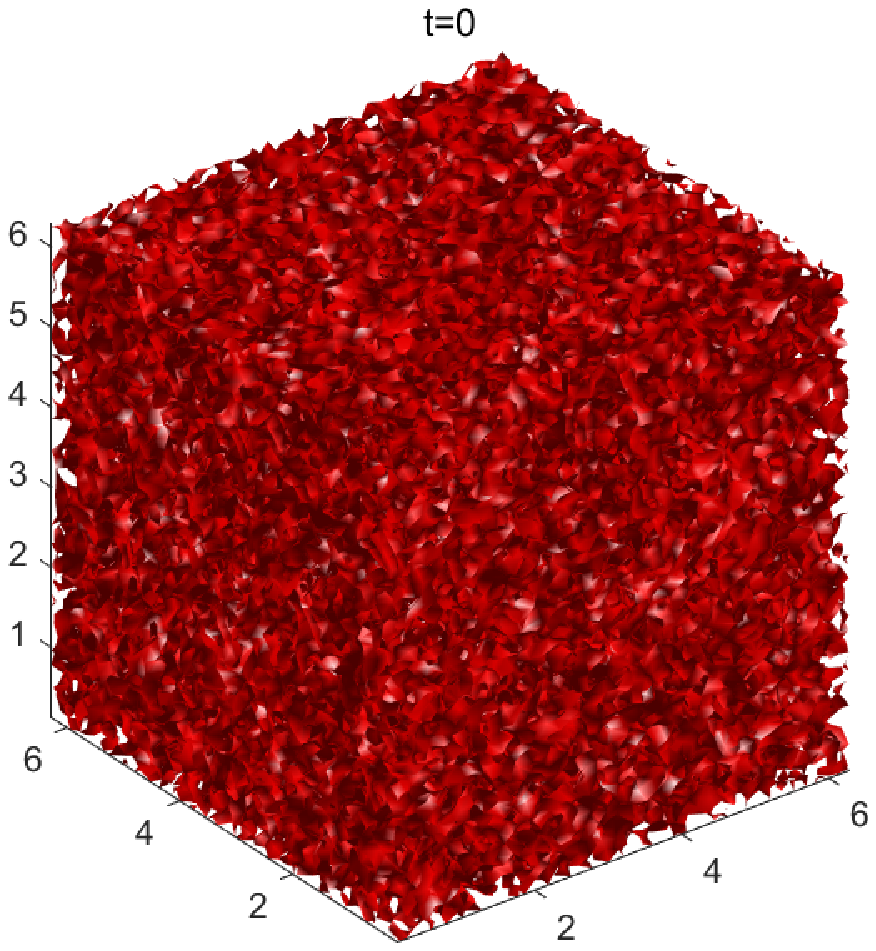}\\
\includegraphics[width=1.3\textwidth]{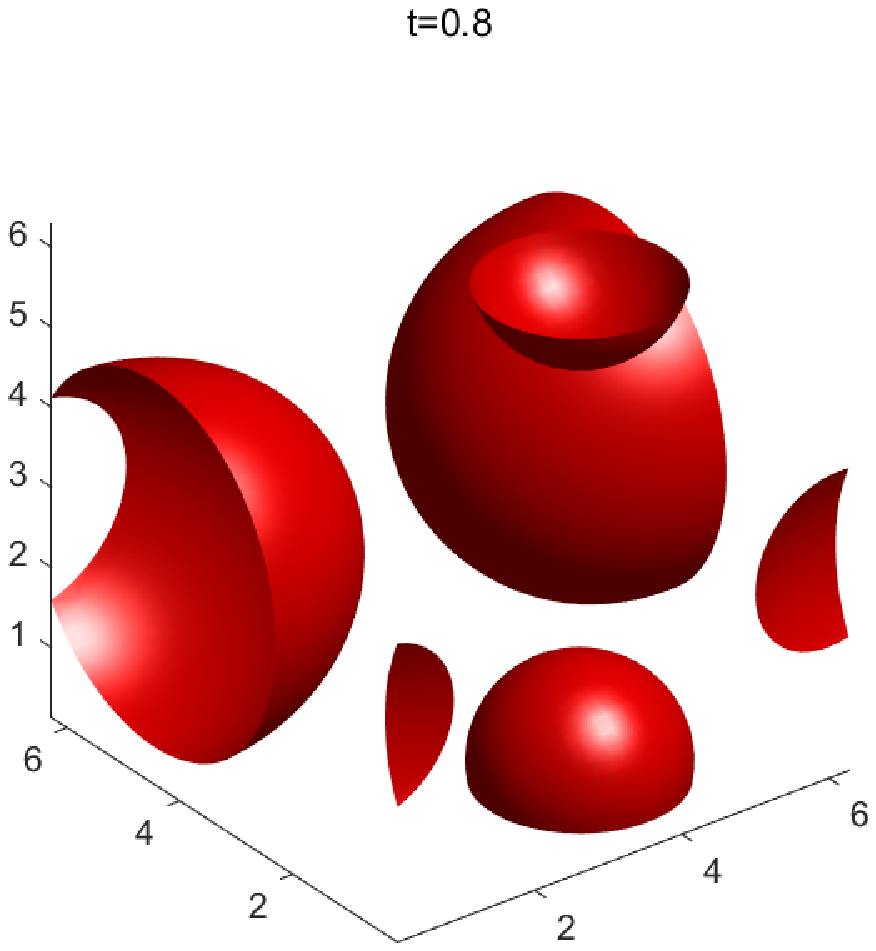}
\end{minipage}
}
\subfigure{
\begin{minipage}[b]{0.31\textwidth}
\includegraphics[width=1.3\textwidth]{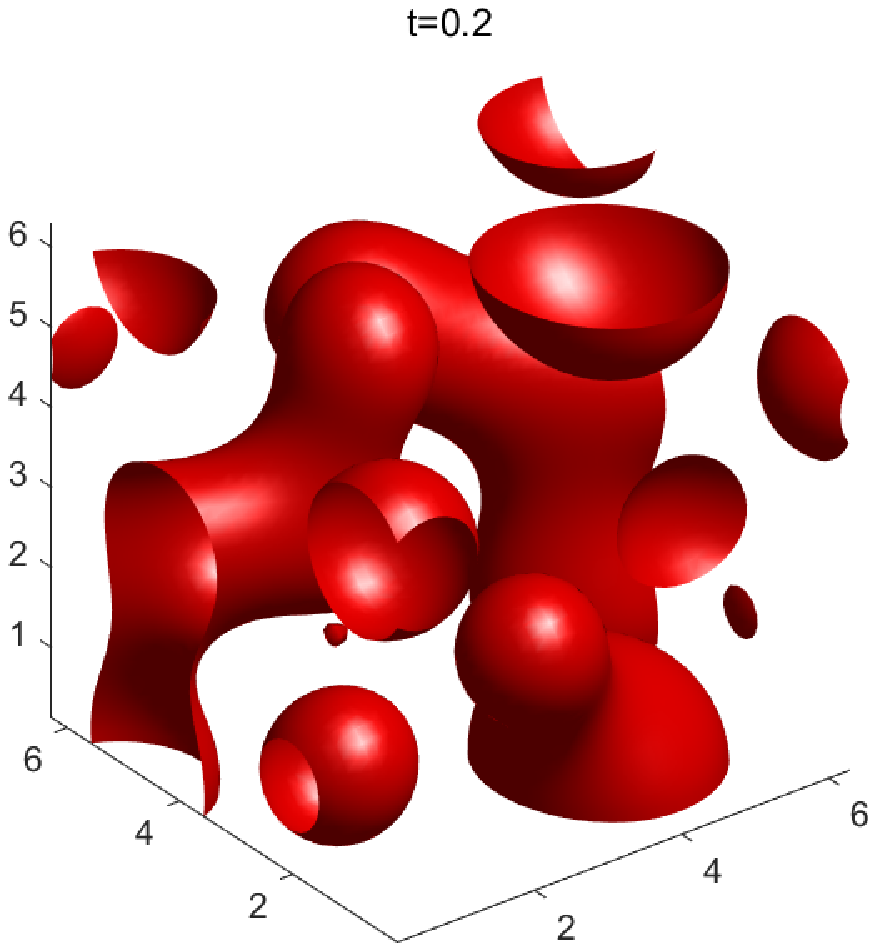} \\
\includegraphics[width=1.3\textwidth]{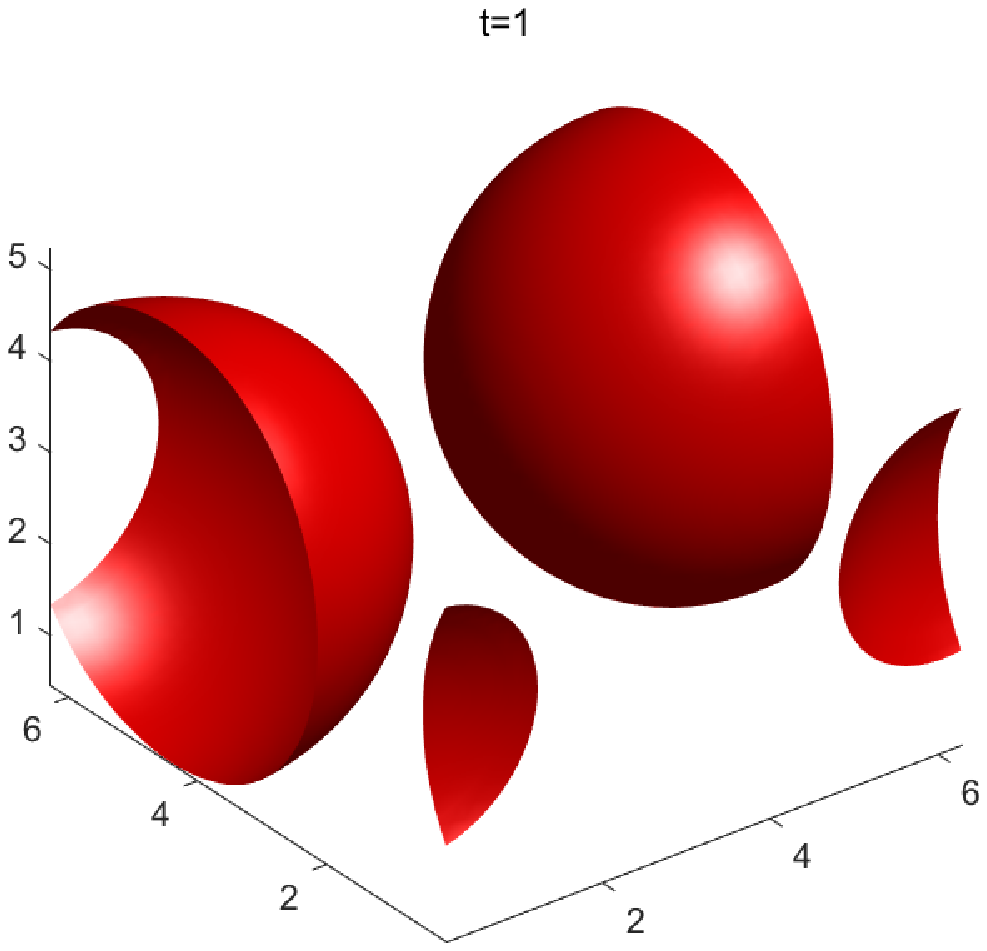}
\end{minipage}
}
\subfigure{
\begin{minipage}[b]{0.31\textwidth}
\includegraphics[width=1.3\textwidth]{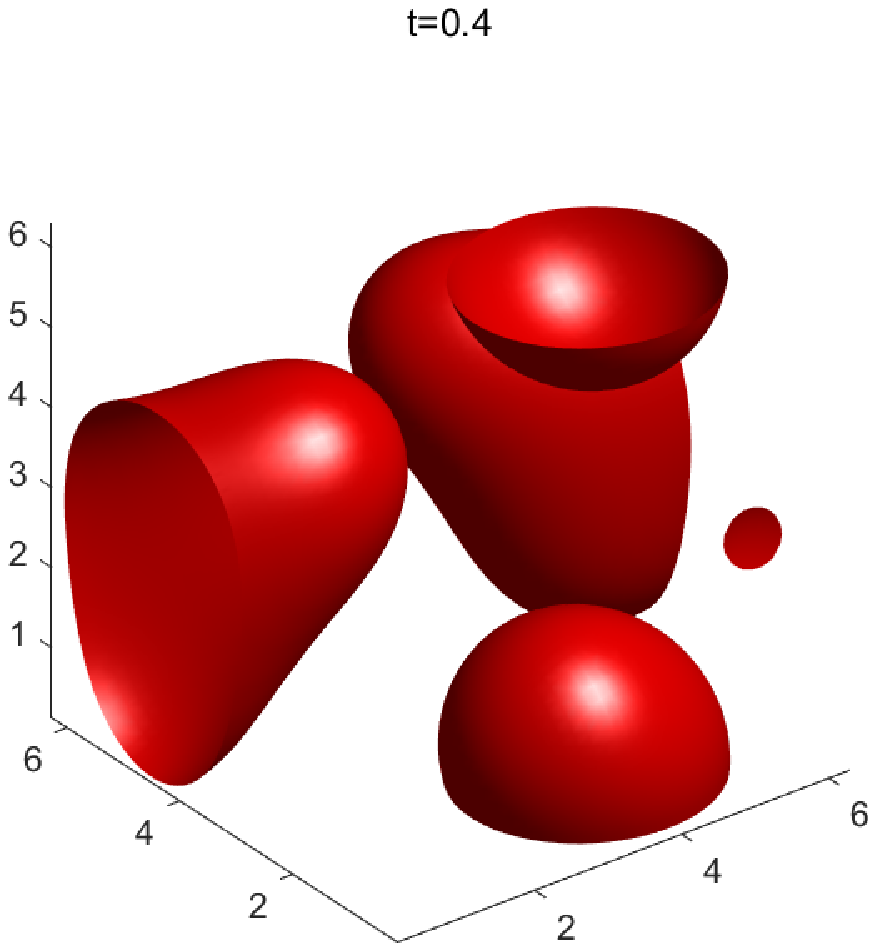} \\
\includegraphics[width=1.3\textwidth]{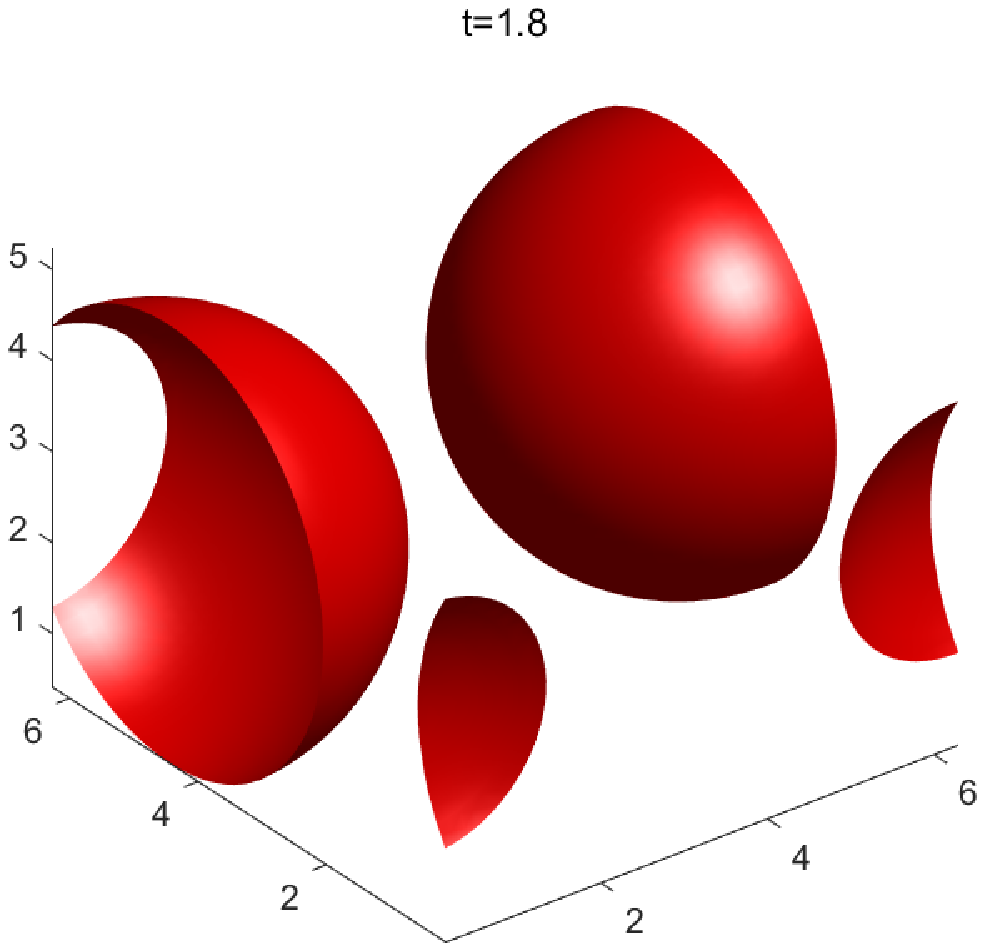}
\end{minipage}
}
\caption{The isosurface of $\{\bm{x} | \Phi(\bm{x})=0\}$ at t = 0,  0.2, 0.4, 0.8, 1, 1.8.}
\label{snap_3d}
\end{figure}

\begin{figure}[!ht]
\centering
\subfigure[]{
\includegraphics[width=6cm]{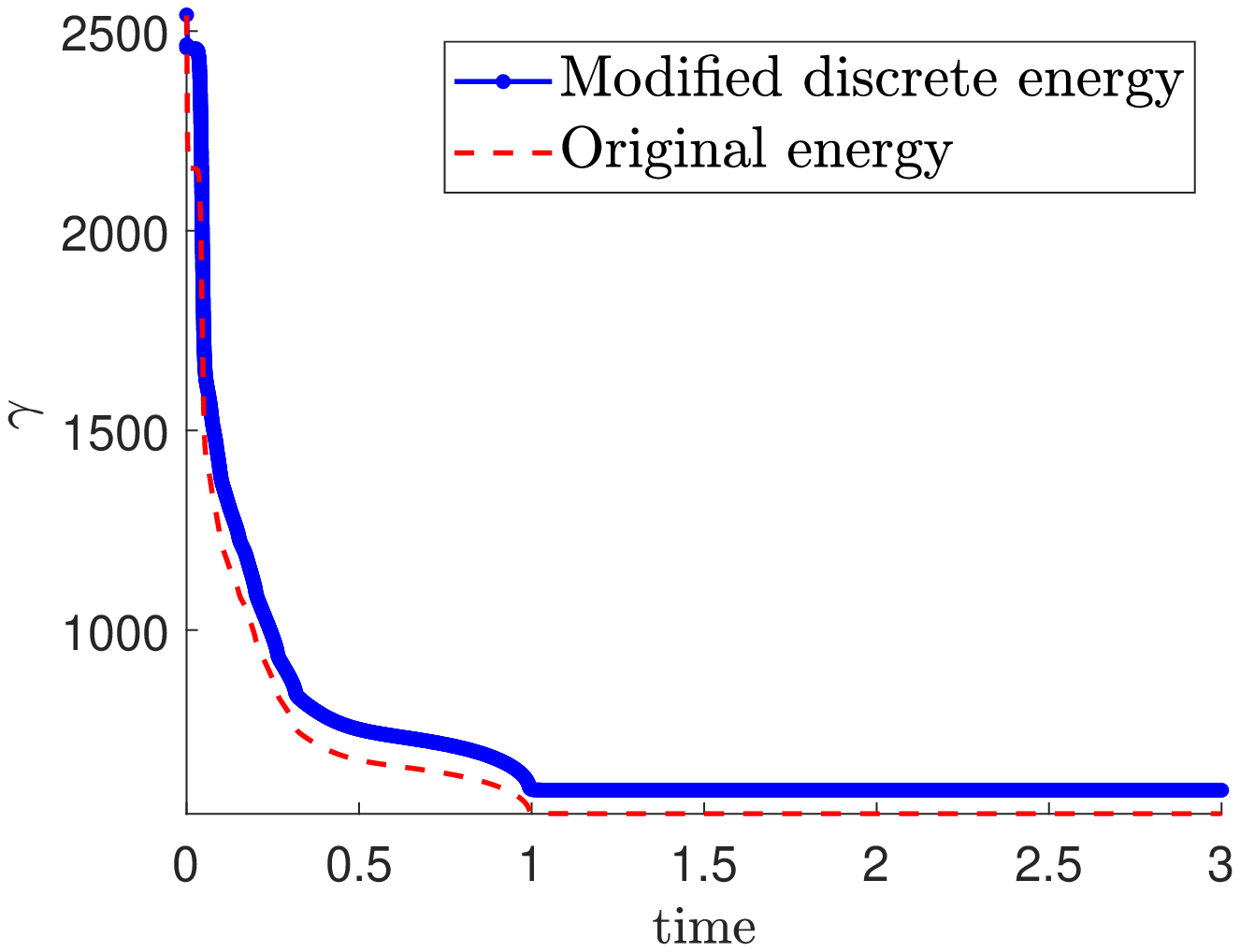}
\label{3D_1}}
\subfigure[]{
\includegraphics[width=6cm]{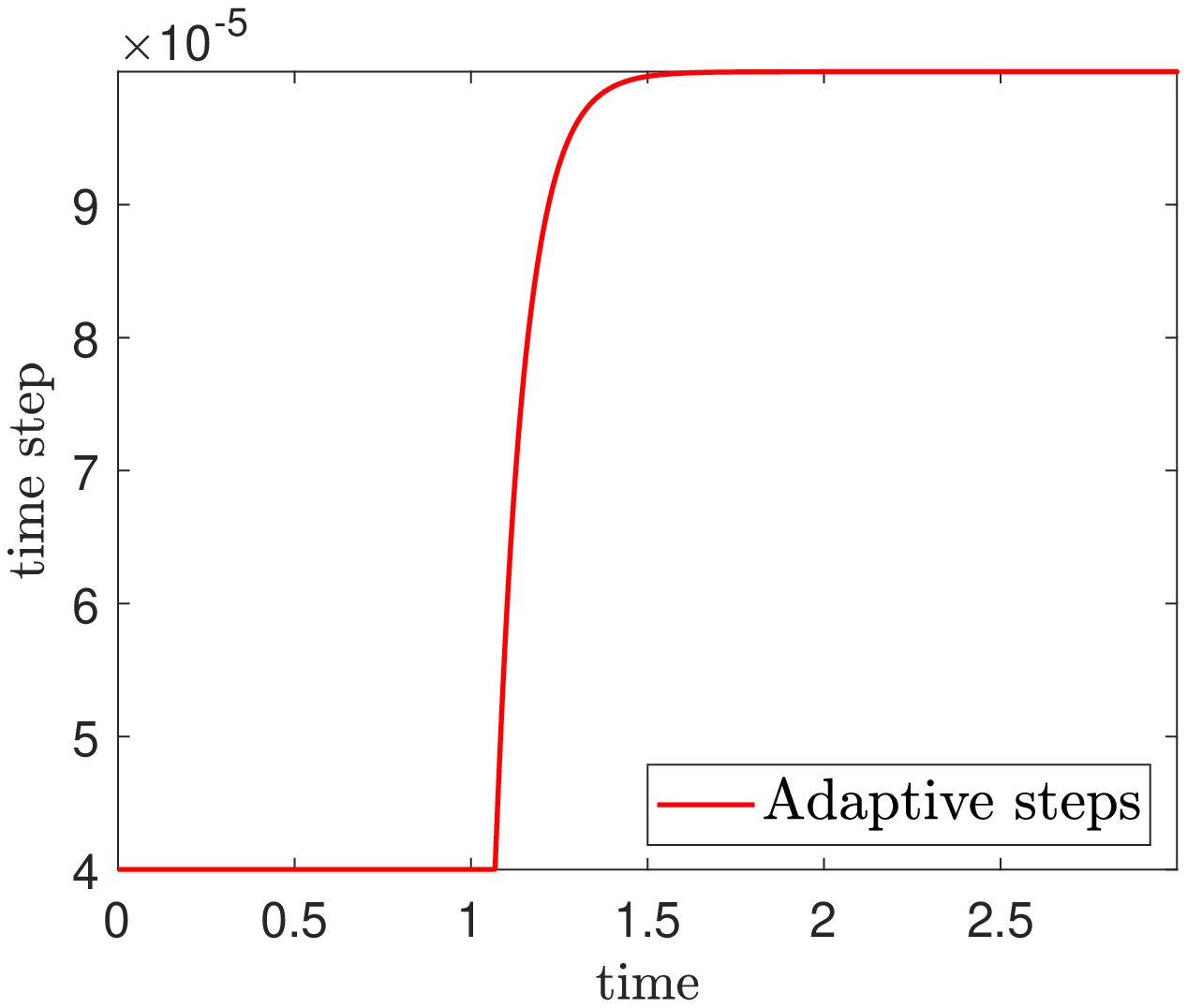}
\label{3D_2}}
\caption{Evolution of modified discrete energy (a) and adaptive time-step size (b) in  3D}
\label{3D}
\end{figure}

The temporal evolution is shown in Figure \ref{snap_3d}. The evolution of phase separation dynamics can be observed. We display the temporal evolutions of original discrete energy and  modified discrete energy in Figure \ref{3D}, which again shows the energy is non-increasing for 3D case.

\begin{figure}[!ht]
\centering
\subfigure[]{
\includegraphics[width=6cm]{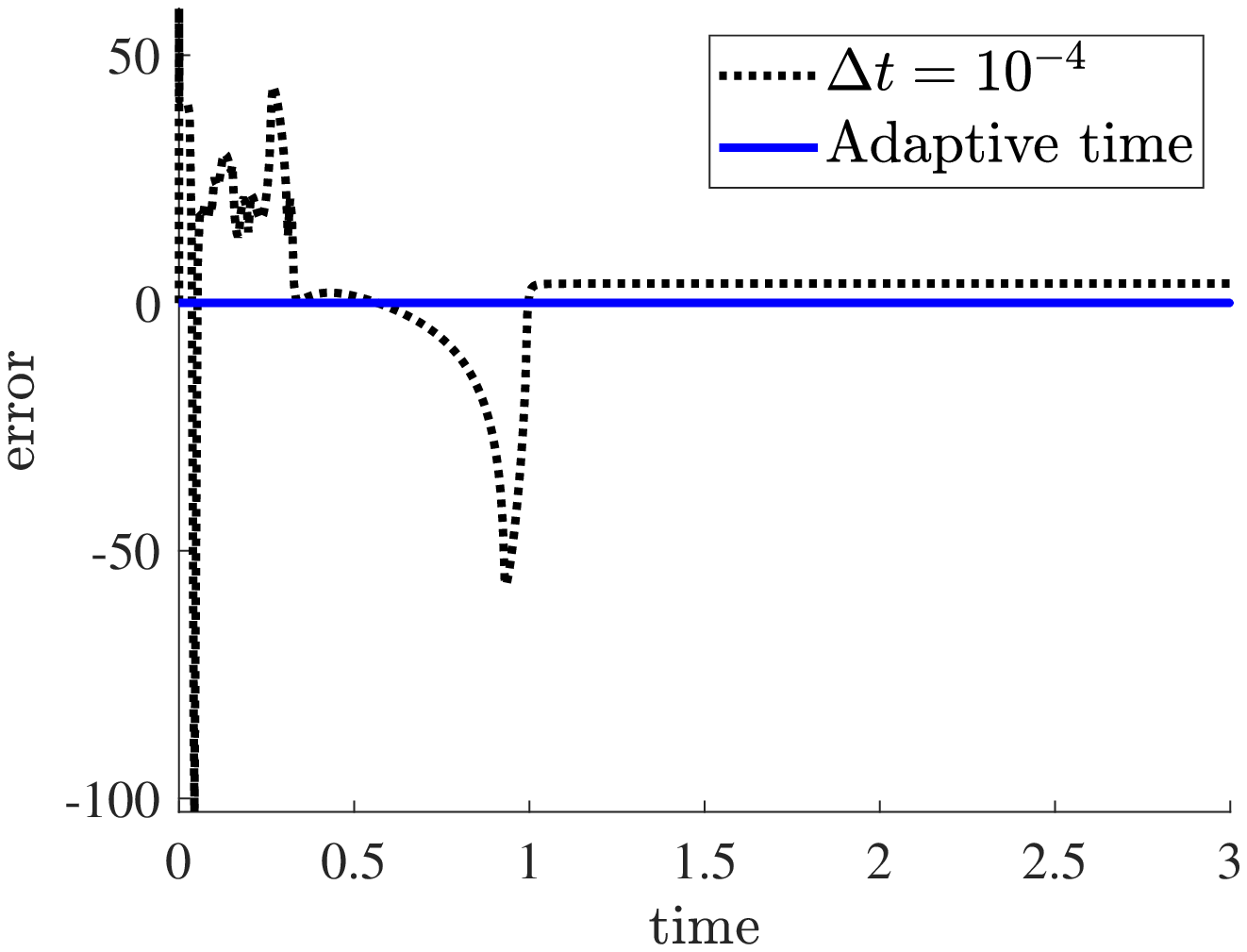}
\label{3D_3}}
\subfigure[]{
\includegraphics[width=6cm]{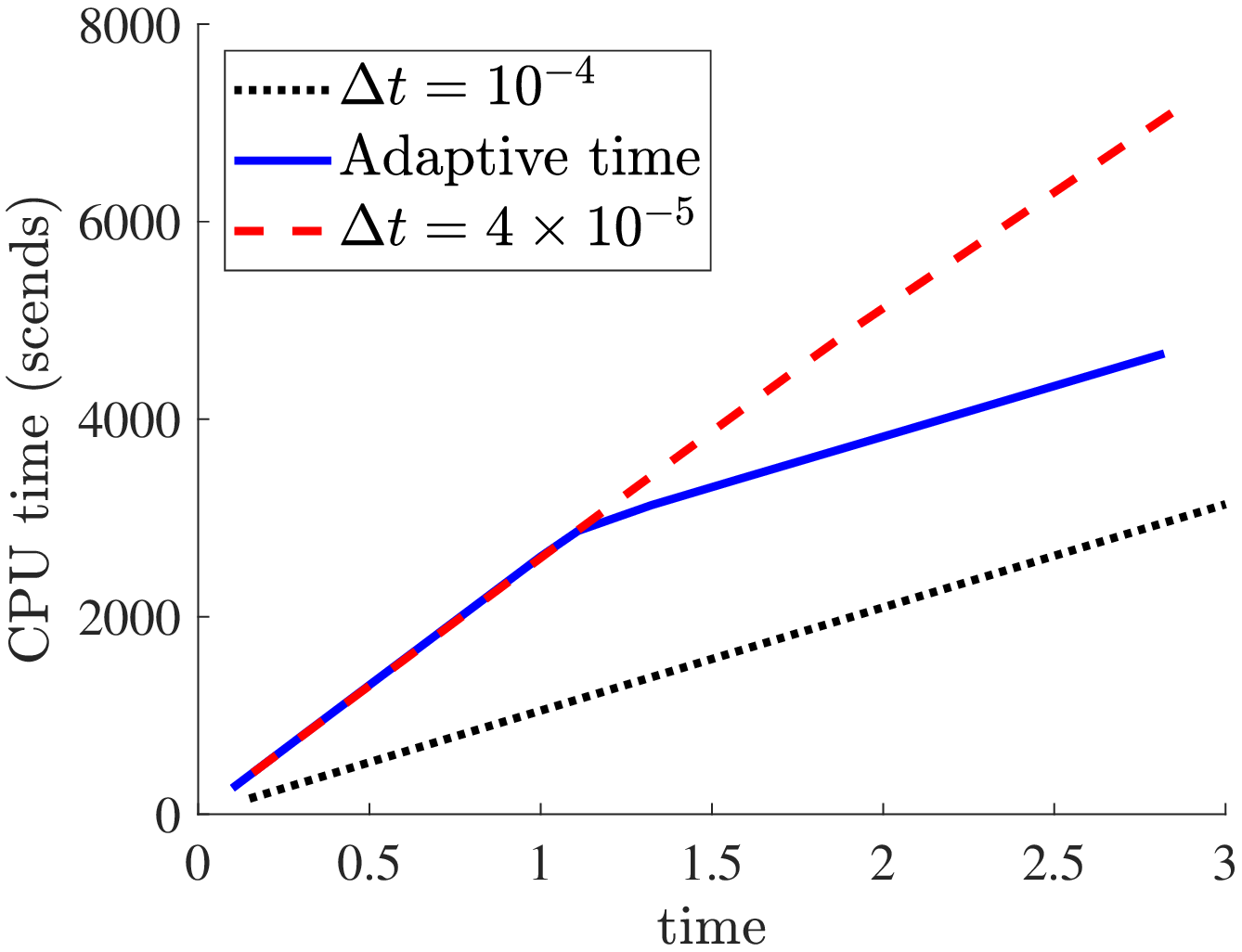}
\label{3D_CPU}}
\caption{Error of $\gamma$ (a) and CPU time contrast (b) in  3D}
\end{figure}

In addition, we study the accuracy and  efficiency of adaptive strategy and fixed-step strategies with time.
Here we take the numerical solution with $\Delta t= 4\times 10^{-5}$ as reference. Figure \ref{3D_3} shows the error evolution of energy of adaptive step size and fixed step size $\Delta t= 1\times 10^{-4}$. In particular, the energy error of the adaptive strategy is always less than $10^{-5}$. Figure \ref{3D_CPU} plots the CPU time, which shows the adaptive strategy significantly reduces CPU time consumption at the almost same accuracy comparing with the very fine fixed time size.

\section{Conclusions}

In this paper, we have constructed and analysed a variable time-step BDF2 IMEX scheme for the C-H equation \eqref{CHeq} by Fourier spectral method in space, which is inspired by the generalized SAV approach introduced in \cite{MR4383075}.  The  unconditional dissipation of the modified energy is proved rigorously.
This paper focuses on the optimal $H^1$-norm error estimation in time under a mild step-ratio condition \Ass{1} (i.e., $0< r_k< 4.8645$).
 To this end, our analysis is mainly based on DOC kernels and its generalized properties. By using the first-order consistent BDF1 scheme to compute the initial step solution $\phi^1$, we achieve the global second-order accuracy of BDF2 scheme with variable time steps. This is obtain by the inequality zoom and some delicate error analysis on the truncation errors.
 We verify numerically that the modified discrete energy of the adaptive scheme is strictly dissipated and related to the original energy.
The proposed adaptive strategy can improve efficiency without sacrificing accuracy, which is verified experimentally by comparing with constant step.
In addition, the obtained results  on variable-step BDF2 can also be extended to some newly developed methods, such as R-GSAV method \cite{MR4097160}, which will be our future work.

\section*{Acknowledgements}
J. Zhang is partially supported by NSFC under grant No. 12171376, 2020-JCJQ-ZD-029, and the Fundamental Research Funds for the Central Universities 2042021kf0050. The numerical simulations in this work have been done on the supercomputing system in the Supercomputing Center of Wuhan University.

\bibliographystyle{abbrv}
\bibliography{add}

\end{document}